\documentclass[3p,times,lefttitle]{elsarticle}

\def\figpath{figures}

\usepackage[T1]{fontenc}
\usepackage[utf8]{inputenc}
\usepackage{amsmath,amssymb,amsthm}
\usepackage{dsfont}
\usepackage{mathtools}
\usepackage{empheq}
\usepackage[binary-units]{siunitx}
\DeclareSIUnit\crate{C}
\usepackage{algorithm}
\usepackage{algpseudocode}
\usepackage{booktabs}
\usepackage{longtable}
\usepackage{doi}
\usepackage[capitalize, noabbrev]{cleveref}
\crefname{subsection}{Subsection}{Subsections}

\usepackage{hyperref}
\hypersetup{
  hidelinks=true
}
\journal{Computers \& Mathematics with Applications}

\biboptions{sort&compress}

\usepackage{amsmath,amssymb}

\makeatletter
\newcommand{\@defs@vec}[1]{\boldsymbol{#1}}
\newcommand{\@defs@tens}[1]{\mathbf{#1}}

\DeclareMathOperator{\RelTol}{RelTol}
\DeclareMathOperator{\AbsTol}{AbsTol}

\newcommand{\de}{\,\mathrm{d}}
\newcommand{\ptl}{\partial}

\newcommand{\IR}{\mathbb{R}}

\newcommand{\vect}[1]{\@defs@vec{#1}{}}
\newcommand{\tens}[1]{\@defs@tens{#1}{}}

\newcommand{\trp}{\textsf{T}}

\newcommand{\grad}{\boldsymbol{\nabla}}

\newcommand{\tfinal}{t_\text{end}}
\newcommand{\sol}{\boldsymbol{y}}

\newcommand{\EH}{\texttt{E}_\text{H}}
\newcommand{\Fo}{\text{Fo}}

\newcommand{\OmegaL}{\Omega_0}

\newcommand{\ch}{\text{ch}}

\newcommand{\el}{\text{el}}
\newcommand{\cmax}{c_{\text{max}}}
\newcommand{\gradL}{\grad_0}
\newcommand{\divgL}{\grad_0 \!\cdot\!}
\newcommand{\something}{\boldsymbol{\cdot}}
\newcommand{\minus}{\scalebox{0.75}[1.0]{$-$}}
\newcommand{\normal}{\vect{n}}
\newcommand{\normalL}{\vect{n}_0}
\newcommand{\displacement}{\vect{u}}
\newcommand{\displacementt}{u}
\newcommand{\displacementD}{\displacement_h}

\newcommand{\possibleDOF}{\mathcal{P}}
\newcommand{\otherDOF}{\mathcal{N}}
\newcommand{\activeset}{\mathcal{A}}
\newcommand{\inactiveset}{\mathcal{I}}
\newcommand{\activesetk}{\activeset_k}
\newcommand{\alldofs}{\mathcal{S}}
\newcommand{\otherandinactive}{{\hat{\mathcal{S}}}}
\newcommand{\Ctensor}{\mathds{C}}

\newcommand{\muD}{\mu_h}
\newcommand{\cD}{c_h}

\newcommand{\mathcalC}{\mathcal{C}}

\newcommand{\FP}[2]
{\frac{\partial #1}{\partial #2}}

\makeatother

\newtheoremstyle{boldremark}
{\dimexpr\topsep/2\relax} % space above
{\dimexpr\topsep/2\relax} % space below
{}          % body font
{}          % indent amount
{\bfseries} % theorem head font
{.}         % punctuation after theorem head
{.5em}      % space after theorem head
{}          % theorem hed spec. (empty = "normal")

\theoremstyle{boldremark}

\usepackage{xcolor}
\usepackage[normalem]{ulem}
\usepackage{cancel}
\newcommand{\stkout}[1]{\ifmmode\text{\sout{\ensuremath{#1}}}\else\sout{#1}\fi}

\begin{document}

  %%%%%%%%%%%%%%%%%%%%%%%%%%%%%%%%%%%%%%%%%%%%%%%%%%%%%%%%%%%%%%%%%%%%%%%%%%%%%%
  %% Frontpage
  %%%%%%%%%%%%%%%%%%%%%%%%%%%%%%%%%%%%%%%%%%%%%%%%%%%%%%%%%%%%%%%%%%%%%%%%%%%%%%
  \begin{frontmatter}

    %% Title
    \title{Simulation of the Deformation for Cycling Chemo-Mechanically Coupled
      Battery Active Particles with Mechanical Constraints}

    %% Authors with affiliation
    \author[ianm]{R. Schoof\corref{correspondingauthor}}
    \cortext[correspondingauthor]{Corresponding author}
    \ead{raphael.schoof@kit.edu}

    \author[tvt]{G. F. Castelli}
    % \ead{fabian.castelli@kit.edu}

    \author[ianm]{W. D\"orfler}
    % \ead{willy.doerfler@kit.edu}

    \address[ianm]{Karlsruhe Institute of Technology~(KIT),
      Institute for Applied and Numerical Mathematics~(IANM),
      Englerstr.~2, 76131~Karlsruhe, Germany}

    \address[tvt]{Karlsruhe Institute of Technology~(KIT),
      Institute of Thermal Process Engineering~(TVT),
      Kaiserstr.~12, 76131~Karlsruhe, Germany}

    %% Abstract
    % !TeX encoding = UTF-8
% !TeX spellcheck = en_US

%%%%%%%%%%%%%%%%%%%%%%%%%%%%%%%%%%%%%%%%%%%%%%%%%%%%%%%%%%%%%%%%%%%%%%%%%%%%%%%%
%% 0-Abstract
%%%%%%%%%%%%%%%%%%%%%%%%%%%%%%%%%%%%%%%%%%%%%%%%%%%%%%%%%%%%%%%%%%%%%%%%%%%%%%%%

% Abstract should be written in the present tense and impersonal style
% (i.e., avoid we), and be at most 200 words long
\begin{abstract}

  %% General
  Next-generation lithium-ion batteries with silicon anodes have positive
  characteristics due to higher energy densities compared to state-of-the-art
  graphite anodes. However, the large volume expansion of silicon anodes
  can cause high mechanical stresses, especially if the battery active particle
  cannot expand freely.
  %% Summary of theroy/numerics
  In this article, a thermodyna\-mi\-cally consistent continuum model for
  coupling
  chemical and mechanical effects of electrode particles is extended by a
  change in the boundary condition for the displacement via a variational
  inequality. This switch represents a limited enlargement of the particle
  swelling or shrinking due to lithium intercalation or deintercalation in the
  host material, respectively. For inequality constraints as boundary condition
  a smaller time step size is need as well as a locally finer mesh. The
  combination of a primal-dual active set algorithm, interpreted as
  semismooth Newton method, and a spatial and temporal adaptive algorithm allows
  the efficient numerical investigation based on a finite element method.
  %% Summary of results
  Using the example of silicon, the chemical and mechanical behavior of one-
  and two-dimensional representative geometries for a charge-discharge cycle is
  investigated. Furthermore, the efficiency of the adaptive algorithm is
  demonstrated. It turns out that the size of the gap has a significant
  influence on the maximal stress value and the slope of the increase.
  Especially in two space dimensions, the obstacle can cause an additional
  region with
  a lithium-poor phase.
\end{abstract}

%%%%%%%%%%%%%%%%%%%%%%%%%%%%%%%%%%%%%%%%%%%%%%%%%%%%%%%%%%%%%%%%%%%%%%%%%%%%%%%%
%% End of 0-abstract.tex
%%%%%%%%%%%%%%%%%%%%%%%%%%%%%%%%%%%%%%%%%%%%%%%%%%%%%%%%%%%%%%%%%%%%%%%%%%%%%%%%

    %% Keywords

  \end{frontmatter}

  %%%%%%%%%%%%%%%%%%%%%%%%%%%%%%%%%%%%%%%%%%%%%%%%%%%%%%%%%%%%%%%%%%%%%%%%%%%%%%
  %% Content
  %%%%%%%%%%%%%%%%%%%%%%%%%%%%%%%%%%%%%%%%%%%%%%%%%%%%%%%%%%%%%%%%%%%%%%%%%%%%%%
  %% Introduction
  % !TeX encoding = UTF-8
% !TeX spellcheck = en_US

%%%%%%%%%%%%%%%%%%%%%%%%%%%%%%%%%%%%%%%%%%%%%%%%%%%%%%%%%%%%%%%%%%%%%%%%%%%%%%%%
%% 1-Introduction
%%%%%%%%%%%%%%%%%%%%%%%%%%%%%%%%%%%%%%%%%%%%%%%%%%%%%%%%%%%%%%%%%%%%%%%%%%%%%%%%

\section{Introduction}
\label{sec:introduction}

%% General introduction
To meet the challenges of climate change, lithium-ion batteries have
emerged as an important and desirable form of energy storage.
The high energy density
and long life time of the electrochemical storage system of lithium-ion
batteries is crucial for mobile applications \cite{tomaszewska2019lithium-ion}.

In addition, batteries with silicon anodes have proven to be very
promising, since their nearly tenfold theoretical capacity compared to
graphite ones currently in
use~\cite{tian2015high, li2021diverting, mo2020tin-graphene}.
However, the additional storage of
lithium-ions can lead to a volume expansion up to~300\%~\cite{zhang2011review}.
The large mechanical stresses occurring as a consequence during the lithiation
and delithiation inside the host material can finally lead to particle fracture
and therefore cause
an undesired shorter battery lifetime and faster aging
process~\cite{xu2016electrochemomechanics, zhao2019review}. Improving the
understanding of the degradation mechanism for lithium-ion batteries with new
materials is an important step towards a sustainable future.

%% Chemical-mechanical coupling
The coupling of chemical and mechanical effects inside the battery active
particles is of great interest to understand the occurrence of the mechanical
and diffusion-induced stress inside the host material \cite{song2015diffusion}.
For example for phase separating materials like lithium manganese oxide
spinel~$\text{Li}_x\text{Mn}_2\text{O}_4$~(LMO), lithium iron
phosphate~$\text{Li}_x\text{FePO}_4$~(LFP) or sodium iron
phosphate~$\text{Na}_x\text{FePO}_4$~(NFP), the stresses are caused by a volume
mismatch between lithium-poor and lithium-rich phases during the intercalation
and deintercalation process~\cite{delmas2008lithium, van-der-ven2000phase,
  walk2014comparison, zhao2019review, song2015diffusion}.

%% Models for large deformation and lithium intercalation
For the coupling of phase separating materials with elastic properties,
the Cahn--Hilliard theory \cite{cahn1958free, cahn1959free} can be extended with
mechanical effects resulting in the Cahn--Larch\'{e}
approach~\cite{larche1973linear, garcke2001cahn-hilliard, garcke2005numerical}
with small deformations and furthermore with finite
deformations~\cite{di-leo2014cahn-hilliard-type, hennessy2020phase,
  walk2014comparison, werner2021multi-field}.
These models have been used in recent years to simulatively investigate the
intercalation of lithium for many different materials, e.g., see for
$\text{Li}_x\text{Mn}_2\text{O}_4$~\cite{walk2014comparison,
  huttin2012phase-field, zhang2018nonlocal},
$\text{Li}_x\text{FePO}_4$~\cite{castelli2021efficient, castelli2021numerical,
  zhang2020mechanically, di-leo2014cahn-hilliard-type, wu2019phase},
$\text{Na}_x\text{FePO}_4$~\cite{zhang2018sodium, zhang2019phase-field,
  zhang2020mechanically} or silicon~\cite{chen2014phase-field,
  zhang2019phase-field_1, poluektov2018modelling,
  kolzenberg2022chemo-mechanical,
  schoof2022parallelization} and the references therein.

%% Change in boundary condition: Primal-dual active set algorithm
All previous simulations of battery active particles have in common that the
considered geometries
can freely swell and are not limited in their volume
enlargement. However, external conditions can change the boundary condition for
the displacement, such as environmental pressure changes or the contact with
the battery case, the current collector or other electrode particles.
This limitation of volume is especially of great significance for the
large volume change of silicon.
In these situations contact problems occur, see for detailed
information~\cite{laursen2002computational,
  wriggers2006computational}.
There are different possibilities to capture such
changes in boundary conditions with various advantages and disadvantages like
\emph{penalty formulation}, \emph{augmented Lagrangian formulation} or
\emph{dual Lagrange multipliers},
compare~\cite[Section~17]{willner2003kontinuums-} and~\cite{alart1991mixed,
  brunssen2007fast, fischer2005frictionless,
  hintermuller2002primal-dual, puso2004mortar, wohlmuth2003monotone}
and the
references therein. The obstacle boundary condition can be written as
\emph{Karush--Kuhn--Tucker}~(KKT) complementary conditions. The equations are
also called \emph{Signorini conditions}, because of their first formulation by
Signorini \cite{signorini1933sopra, signorini1933sopra_1} for the unilateral
normal contact.

The \emph{primal-dual active set strategy} is an efficient technique for
this kind of inequality constraints
and can be interpreted as
\emph{semismooth Newton method}, compare \cite{hueber2005primal-dual,
  hueber2005priori, hueber2013contact, hintermuller2002primal-dual,
  hintermuller2003semismooth, frohne2016efficient, hager2010semismooth} and
\texttt{deal.II} tutorial step-41 in \cite{arndt2021deal-ii}. An additional
ansatz with a direct approach for Signorini's problem with small deformations
and linear elasticity is investigated in \cite{kornhuber2001adaptive}. In case
of time-dependent problems coupled with inequality boundary constraints there
are several numerical solution approaches \cite{hager2010semismooth,
  sander2013towards, de-los-reyes2012combined, lauser2011new,
  sa-ngiamsunthorn2021optimal}. The dual Lagrange multiplier ansatz has the
advantage that there is no need to change the system size. Furthermore, this
method fulfills the obstacle boundary inequality constraints in the weak
integral sense and the condition number of the system matrix does not
change~\cite{brunssen2007fast}.

%% Numerical challenges
The numerical simulation of battery active particles with limited swelling
due to inequality constrains at the boundary is computationally challenging
because the contact region, which requires a higher grid resolution, changes
in time.
Space and time adaptivity is crucial to properly capture all
relevant effects. Additionally, the switch from charging to discharging a
lithium-ion battery for long term cycle investigations needs
an appropriate mechanism for space and time control.

%% Review on used methods in literature
In \cite{chen2014phase-field, zhang2019phase-field_1}, a phase separation
ansatz is used to model the intercalation of silicon. Following
\cite{di-leo2015diffusion-deformation}, however, a two-phase lithiation
mechanism occurs only in the first half cycle of lithiation and therefore will
not be considered in this paper. In particular, in
\cite{kolzenberg2022chemo-mechanical} a measured open circuit voltage~(OCV)
curve is used for the chemical energy density in combination with a linear
elastic approach to model the elastic deformation. A further possibility
would be a
Neo-Hookean ansatz as in \cite{werner2021multi-field} for the
mechanical deformation.

%% Our findings
In this article, we rely on the developed model approach by
\cite{kolzenberg2022chemo-mechanical} for silicon with a finite deformation
ansatz. We combine the large volume expansion with the inequality constraints
for an obstacle boundary problem to
simulate particle swelling in a limited space.
The primal-dual active set strategy is derived, applied as semismooth
Newton method and added to the spatial and temporal adaptive solution algorithm
by \cite{castelli2021numerical, castelli2021efficient, castelli2021study} with
higher finite element order and a fully variable order, variable time step size
time integration scheme. This makes it possible to consider various parameter
setups and one- and two-dimensional
geometrical setups due to large computational
savings introduced by the adaptive algorithm.

%% Organization of this article
The remaining paper is structured as follows: in \cref{sec:theory}, we present
our model to characterize the chemical-mechanical coupling with the obstacle
boundary condition during one charging and discharging cycle. Next, we derive
the semismooth Newton method from the primal-dual active set algorithm
and combine it with a space and time adaptive algorithm. \cref{sec:results}
discusses the simulation results for the developed model with obstacle contact.
Finally, we summarize our main findings in \cref{sec:conclusion}.

%%%%%%%%%%%%%%%%%%%%%%%%%%%%%%%%%%%%%%%%%%%%%%%%%%%%%%%%%%%%%%%%%%%%%%%%%%%%%%%%
%% End of 1-introduction.tex
%%%%%%%%%%%%%%%%%%%%%%%%%%%%%%%%%%%%%%%%%%%%%%%%%%%%%%%%%%%%%%%%%%%%%%%%%%%%%%%%

  %% Modeling
  % !TeX encoding = UTF-8
% !TeX spellcheck = en_US

%%%%%%%%%%%%%%%%%%%%%%%%%%%%%%%%%%%%%%%%%%%%%%%%%%%%%%%%%%%%%%%%%%%%%%%%%%%%%%%%
%% 2-Theory
%%%%%%%%%%%%%%%%%%%%%%%%%%%%%%%%%%%%%%%%%%%%%%%%%%%%%%%%%%%%%%%%%%%%%%%%%%%%%%%%

%%%%%%%%%%%%%%%%%%%%%%%%%%%%%%%%%%%%%%%%%%%%%%%%%%%%%%%%%%%%%%%%%%%%%%%%%%%%%%%%
\section{Theory}
\label{sec:theory}
%%%%%%%%%%%%%%%%%%%%%%%%%%%%%%%%%%%%%%%%%%%%%%%%%%%%%%%%%%%%%%%%%%%%%%%%%%%%%%%%

In this section, we review and summarize the theory from
\cite{kolzenberg2022chemo-mechanical, castelli2021efficient, brunssen2007fast}
to formulate the coupled chemical and mechanical particle obstacle problem for
battery active particles. For this we state a thermodynamically consistent
theory for the chemo-mechanical coupling for (de-)lithiation with inequality
boundary constraints to an obstacle problem. In a first step, we introduce the
finite deformation theory for the particle and then couple the chemical and
mechanical effects with a common free energy density. After derivation of the
equations for chemistry and mechanics, we incorporate the boundary constraints
for the representation of an obstacle hindering the particle to expand freely.
Since we consider the intercalation and deintercalation of lithium into and out
of the host material, we simplify our wording and combine \emph{lithiation and
  delithiation} as well as \emph{charging and discharging} in the word
\emph{cycling}. A selection of abbreviations and symbols of our work is listed
in~\cref{app:abbreviations_and_symbols} and some notation explanations notation
are given in~\cref{app:tensor_analysis}.

%%%%%%%%%%%%%%%%%%%%%%%%%%%%%%%%%%%%%%%%%%%%%%%%%%%%%%%%%%%%%%%%%%%%%%%%%%%%%%%%
\subsection{Finite Deformation}
\label{subsec:finite_deformation}
%%%%%%%%%%%%%%%%%%%%%%%%%%%%%%%%%%%%%%%%%%%%%%%%%%%%%%%%%%%%%%%%%%%%%%%%%%%%%%%%

\begin{figure}[!bt]
  \centering
  \includegraphics[width = 0.55\textwidth,
  page=1]{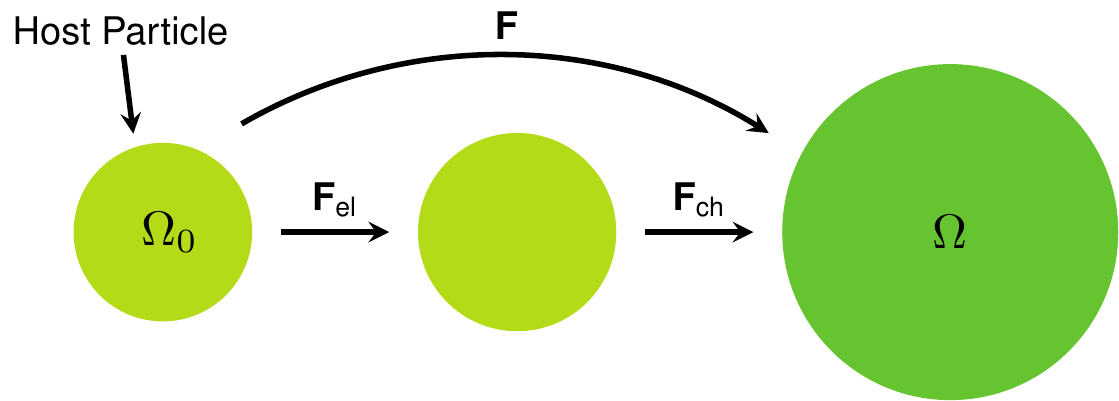}
  \caption{The total deformation $\tens{F}$ can be
    multiplicatively decomposed into an elastic part $\tens{F}_\el$ and a
    chemical part~$\tens{F}_\ch$,
    compare~\cite[Figure~1]{schoof2022parallelization}.}
  \label{fig:deformation_theory}
\end{figure}

To model the particle deformation during cycling, we consider a
motion~$\vect{x}(t, \vect{X}_0) =\vect{\chi}(t, \vect{X}_0) =
\vect{X}_0 + \vect{u}(t, \vect{X}_0)$
with a mapping
$\vect{\chi}\colon \IR_{\geq 0} \times \Omega_0 \to \Omega$. Here, $\vect{X}_0
\in \Omega_0$ corresponds to an arbitrary point in the Lagrangian reference
configuration~$\Omega_0$ which is mapped to a point $\vect{x} \in \Omega$ in
the current Eulerian configuration~$\Omega$. The reversible total deformation
gradient tensor $\tens{F}$ is defined as $\tens{F}= \ptl \vect{\chi}/\ptl
\vect{X_0}$ \cite[Chapter 2.4]{holzapfel2000nonlinear}. This results in the
relation $\tens{F} = \tens{Id} + \gradL \vect{u}$ with the identity
matrix~$\tens{Id}$ and the gradient of the displacement $\gradL \vect{u}$ with
respect to the spatial coordinates of the reference configuration,
compare~\cite[Chapter~VI\S{}1]{braess2007finite} and
\cite[Section~2.4]{holzapfel2000nonlinear}. We follow
\cite{kolzenberg2022chemo-mechanical} and multiplicatively decompose the
deformation gradient as~$\tens{F} = \tens{F}_\ch \tens{F}_\el$.
A sketch of this decomposition is given in \cref{fig:deformation_theory}. The
elastic part $\tens{F}_\el$ occurs due to mechanical stress, whereas the
chemical part results from the changes in the lithium concentration. With an
isotropic and linear chemical expansion of the active material the chemical
part of the deformation gradient is given by $\tens{F}_\ch = \lambda_\ch
\tens{Id}$ with $\lambda_\ch = \sqrt[3]{1+v_\text{pmv}c}$, where $v_\text{pmv}$
defines the partial molar volume of the host material and $c$ the lithium
concentration~\cite{castelli2021efficient}.

%%%%%%%%%%%%%%%%%%%%%%%%%%%%%%%%%%%%%%%%%%%%%%%%%%%%%%%%%%%%%%%%%%%%%%%%%%%%%%%%
\subsection{Free Energy}
\label{subsec:free_energy}
%%%%%%%%%%%%%%%%%%%%%%%%%%%%%%%%%%%%%%%%%%%%%%%%%%%%%%%%%%%%%%%%%%%%%%%%%%%%%%%%
Based on a free energy density $\psi$, we use a thermodynamically consistent
model to guarantee a strictly positive entropy
production~\cite{latz2015multiscale, latz2011thermodynamic,
  kolzenberg2022chemo-mechanical, schammer2021theory}. Following
\cite{kolzenberg2022chemo-mechanical}, we define the free energy density~$\psi$
as
\begin{align}
  \label{eq:free_energy_density}
  \psi(c, \grad_0 \vect{u}) = \psi_\ch(c) + \psi_\el(c, \grad_0 \vect{u}),
\end{align}
combining chemical and mechanical effects. Adding an interfacial part
$\psi_\text{int}(\grad c)$ to \cref{eq:free_energy_density} for materials with
phase separation leads to the classical Cahn--Hilliard approach combined with
elasticity \cite{castelli2021efficient, anand2012cahn-hilliard-type,
  di-leo2014cahn-hilliard-type, zhang2018lithiation-induced, zhang2018sodium}.

For the definition of the chemical part $\psi_\ch$ we use the experimentally
obtained OCV curve~$U_\text{OCV}$ \cite{chan2007high-performance,
  kolzenberg2022chemo-mechanical, keil2016calendar, latz2015multiscale,
  latz2013thermodynamic}
\begin{align}
  \psi_\ch(c) = \minus\int_0^{c / c_{\max}} F \, U_\text{OCV}(z) \de z
\end{align}
with the Faraday constant $F$. For the elastic part $\psi_\el$ we use the
linear elastic approach (Saint Venant--Kirchhoff model) as in \cite[Section
6.5]{holzapfel2000nonlinear}, \cite[Chapter VI \S{}3]{braess2007finite} and
\cite{castelli2021efficient, kolzenberg2022chemo-mechanical}
\begin{align}
  \psi_\el = \frac{1}{2}
  \tens{E}_\el\!:\!\Ctensor \left[ \tens{E}_\el \right] \qquad
  \text{with} \qquad \Ctensor \left[ \tens{E}_\el \right] = \lambda_\text{H}
  \text{tr}(\tens{E}_\el)\tens{Id} + 2 G_\text{H} \tens{E}_\el,
\end{align}
first and second Lam\'{e} constants $\lambda_\text{H} = 2 G_\text{H} \nu
/\left(1-2\nu\right)$ and $G_\text{H} = \EH/\left[2
\left(1+2\nu\right)\right]$, Young's modulus $\EH$ and Poisson's ratio~$\nu$.
Furthermore, we define the elastic strain tensor~$\tens{E}_\el$, also called
Green--Lagrange strain tensor, in our model by
\begin{align}
  \tens{E}_\el = \frac{1}{2} \left(\tens{F}_\el^{\trp}\tens{F}_\el -
  \tens{Id}\right) =
  \frac{1}{2} \left(\lambda_\ch^{\minus 2}\tens{F}^{\trp}\tens{F} -
  \tens{Id}\right).
\end{align}

%%%%%%%%%%%%%%%%%%%%%%%%%%%%%%%%%%%%%%%%%%%%%%%%%%%%%%%%%%%%%%%%%%%%%%%%%%%%%%%%
\subsection{Elastic Deformation}
%%%%%%%%%%%%%%%%%%%%%%%%%%%%%%%%%%%%%%%%%%%%%%%%%%%%%%%%%%%%%%%%%%%%%%%%%%%%%%%%

A momentum balance law governs the mechanical behavior for the deformation in
the Lagrangian frame
\begin{align}
  \vect{0} &= \minus \divgL \tens{P} \qquad \text{in } (0, \tfinal) \times
  \OmegaL
\end{align}
without considering any body or inertial forces \cite{castelli2021efficient,
  kolzenberg2022chemo-mechanical}. The first Piola--Kirchhoff stress
tensor~$\tens{P}$ and the Cauchy stress $\boldsymbol{\sigma}$ in the Eulerian
frame are coupled via $\tens{P} =
\det{\left(\tens{F}\right)\boldsymbol{\sigma}\tens{F}^{\minus \trp}}$
\cite[Section~3.1]{holzapfel2000nonlinear} using Nanson's formula for a vector
element of infinitesimally small surface area
\cite[Section~2.4]{holzapfel2000nonlinear}. A thermodynamically consistent
derivation specifies the first Piola--Kirchhoff stress
tensor~$\tens{P}(c, \grad_0 \vect{u}) = \partial_\tens{F}\psi=
\lambda_\ch^{\minus 2} \tens{F} \Ctensor \left[ \tens{E}_\el \right]$,
compare~\cite[Section~6.1]{holzapfel2000nonlinear}.

%%%%%%%%%%%%%%%%%%%%%%%%%%%%%%%%%%%%%%%%%%%%%%%%%%%%%%%%%%%%%%%%%%%%%%%%%%%%%%%%
\subsection{Chemistry}
\label{subsec:chemistry}
%%%%%%%%%%%%%%%%%%%%%%%%%%%%%%%%%%%%%%%%%%%%%%%%%%%%%%%%%%%%%%%%%%%%%%%%%%%%%%%%

A continuity equation is used to describe the change of the lithium
concentration inside the host material via
\begin{align}
  \partial_t c = \minus \divgL\boldsymbol{N} \qquad \text{in } (0, \tfinal)
  \times
  \OmegaL
\end{align}
with the lithium flux $\vect{N} \coloneqq \minus {m}(c, \grad_0\vect{u})
\grad_0{\mu} = \minus D \, \bigl( \ptl_c {\mu} \bigr)^{-1} \grad_0{\mu}$ and
the diffusion coefficient $D$ for lithium atoms inside the active material. The
definition for the lithium flux $\vect{N}$
follows~\cite{kolzenberg2022chemo-mechanical} to guarantee positive entropy
production. The chemical potential $\mu = \partial_c \psi$ is stated as the
variational derivative of the Ginzburg--Landau free energy~$\Psi =
\int_{\Omega_0} \psi \de \vect{X}_0$~\cite{hoffmann2018influence}. This leads
to the definition of the chemical potential
\begin{align}
  \mu = \ptl_c \psi &= \minus F U_{\text{OCV}} -
  \frac{v_\text{pmv}}{3}\lambda_\ch^{\minus 5} \left(\tens{F}^\trp
  \tens{F}\right)
  \!:\! \Ctensor \left[ \tens{E}_\el \right]
  = \minus F U_{\text{OCV}} - \frac{v_\text{pmv}}{3 \lambda_\ch^3} \tens{P}
  \!:\!\tens{F}.
\end{align}
The representative particle is cycled with a uniform and constant external flux
$N_\text{ext}$ with either positive or negative sign. This external flux is
applied at the boundary of $\Omega_0$ and measured in terms of the C-rate, for
which we refer to~\cite{castelli2021efficient}. With this definition, the
simulation time $t$ and the state of charge (SOC) can be related by
\begin{align}
  \mathrm{SOC} = \frac{1}{V_{\OmegaL}}\int_{\Omega_0} \frac{c}{c_\text{max}}
  \de \vect{X}_0
  = \frac{c_0}{c_\text{max}} + N_\text{ext} \left[\si{C}\right] \cdot t
  [\si{\hour}]
\end{align}
with the volume $V_{\Omega_0}$ of $\Omega_0$ and a constant initial condition
$c_0 \in \left[0, c_{\max}\right]$.

%%%%%%%%%%%%%%%%%%%%%%%%%%%%%%%%%%%%%%%%%%%%%%%%%%%%%%%%%%%%%%%%%%%%%%%%%%%%%%%%
\subsection{Obstacle Contact Problem}
%%%%%%%%%%%%%%%%%%%%%%%%%%%%%%%%%%%%%%%%%%%%%%%%%%%%%%%%%%%%%%%%%%%%%%%%%%%%%%%%
In the situation of a freely expanding particle during cycling, a stress-free
boundary condition in normal direction is assumed on the particle surface
\cite{castelli2021efficient, schoof2022parallelization}:
\begin{align}
  \label{eq:stress-free_boundary_condition}
  \tens{P} \cdot \vect{n}_0 = \vect{0} \qquad \text{on } (0, \tfinal) \times
  \ptl
  \Omega_0,
\end{align}
where $\vect{n}_0$ is the outer unit vector of the reference
configuration~$\Omega_0$. However, if we consider a situation where the
particle can no longer expand freely, this boundary condition must be adapted.
This means that we have to handle a contact problem as restriction for the
particle swelling.

\begin{figure}[!tb]
  \centering
  \includegraphics[width = 0.8\textwidth,
  page=1]{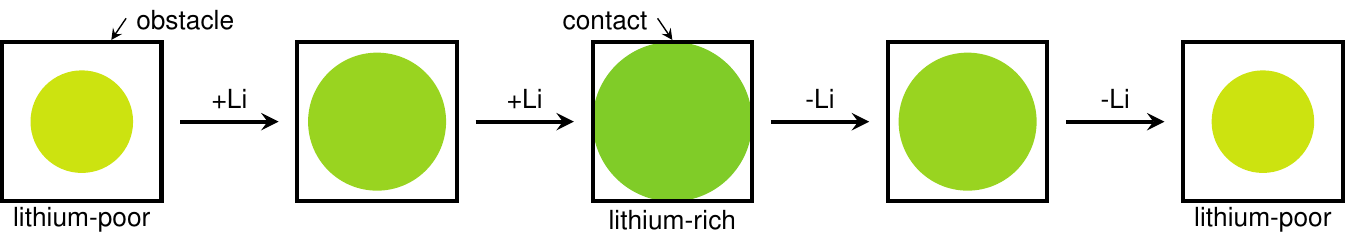}
  \caption{Schematic sketch of lithiation and following delithiation of a
    representative battery active particle with volume change, getting in
    contact with the obstacle and detaching from the obstacle again.}
  \label{fig:schematic-contact-lisi_lithiation_delithiation}
\end{figure}

A schematic sketch of a lithiation and delithiation cycle is shown in
\cref{fig:schematic-contact-lisi_lithiation_delithiation}: for example, the
cross section of a lithium-poor particle is surrounded by a square shaped
obstacle.
The particle increases until it gets in contact with the obstacle.
Now, the stress-free boundary condition in normal direction is replaced by a
restriction of the displacement and thus, nonzero stresses in normal direction
are possible. During delithiation, the particle detaches from the obstacle and
shrinks again until it returns to a lithium-poor state.
In the following, we replace the stress-free boundary condition in
\cref{eq:stress-free_boundary_condition} by an appropriate condition to
incorporate the obstacle contact. As in \cite{hintermuller2003semismooth,
  hintermuller2002primal-dual} we take the new boundary condition:
\begin{subequations}
  \label{eq:contact_boundary_omega}
  \begin{empheq}[left=\empheqlbrace]{alignat=4}
    \displacement &\leq \vect{g}
    \qquad && \text{on }
    (0, \tfinal) \times \partial\Omega, \label{eq:contact_boundary_omega_g}\\
    \boldsymbol{\sigma} \cdot \normal &\leq \vect{0} && \text{on }
    (0, \tfinal) \times \partial\Omega,
    \label{eq:contact_boundary_omega_stress}\\
    \left[\boldsymbol{\sigma} \cdot \normal\right]
    \left[\displacement- \vect{g} \right]
    &= \vect{0} && \text{on }
    (0, \tfinal) \times \partial\Omega,
    \label{eq:contact_boundary_omega_complementary}
  \end{empheq}
\end{subequations}
understood componentwise as introduced in~\cref{app:tensor_analysis}, with
outer unit vector $\vect{n}$ of the current configuration $\Omega$.

Let~$\hat{G}$ be the set that encloses $\Omega_0 $, i.e., $\Omega_0 \subset
\hat{G}$, with boundary $\partial\hat{G}$ that defines the (time-independent)
obstacle. We define
$\hat{\vect{g}}\colon \partial \Omega_0 \rightarrow \partial \hat{G}$,
$\vect{X}_0
\mapsto \hat{\vect{g}}\left(\vect{X}_0\right)$ as the
projection of an arbitrary point~$\vect{X}_0$ on $\partial \Omega_0$ to
the nearest boundary point on~$\partial \hat{G}$ for each component parallel
to the coordinate axes.
Now, we set the \textit{gap function} to
$\vect{g} \colon \partial \Omega_0 \rightarrow \IR_{>
  0}^d$,
$\vect{g}\left(\vect{X}_0\right) \coloneqq \hat{\vect{g}}(\vect{X}_0) -
\vect{X}_0$,
indicated in~\cref{fig:particle_obstacle_2d}.
\cref{eq:contact_boundary_omega_g} thus means that the time dependent
displacement~$\vect{u}$
must componentwise
be smaller or equal to the gap function~$\vect{g}$ between
the
particle and the obstacle.
Furthermore, \cref{eq:contact_boundary_omega_stress} expresses the fact that
for each component, the Cauchy stress tensor in
normal direction is zero or less than zero. The latter case indicates
compressive stress for the respective coordinate component.
The last \cref{eq:contact_boundary_omega_complementary} is
the complementary condition and specifies that one of the two
\crefrange{eq:contact_boundary_omega_g}{eq:contact_boundary_omega_stress} must
be zero, whereas the other condition can be nonzero.
Again, this is interpreted componentwise for each component of the
underlying coordinate system.
For example, it is also possible for the complementary
condition, that in the first component the displacement equals the gap
function and there
is compressive stress in the first component of the Cauchy stress
vector in normal direction, but in the second component the
Cauchy stress vector in normal direction is zero and the displacement
component is smaller than the gap function component.
In short: if the particle
is not in contact with the obstacle, there must be zero stress in
each component of the Cauchy stress tensor in normal
direction, or vice versa, if there are compressive stresses, the particle has
to be in contact with the obstacle in each component.
This kind of boundary obstacle problem is
also known as \emph{Signorini problem} or {thin obstacle
  problem}~\cite[Section~1.11]{friedman1982variational}.

Since we formulate all constitutive equations in the Lagrangian domain,
we use Nanson's formula and it follows
for \cref{eq:contact_boundary_omega}:
\begin{subequations}
  \label{eq:contact_boundary_omega_lagrangian}
  \begin{empheq}[left=\empheqlbrace]{alignat=4}
    \displacement - \vect{g} &\leq \vect{0}
    \qquad
    &&
    \text{on } (0, \tfinal)
    \times
    \partial\Omega_0,
    \label{eq:contact_boundary_omega_g_lagrangian}\\
    \minus \tens{P} \cdot \normal_0 &\geq \vect{0} && \text{on } (0, \tfinal)
    \times
    \partial\Omega_0,
    \label{eq:contact_boundary_omega_stress_lagrangian}\\
    \left[\minus \tens{P} \cdot \normal_0\right]
    \left[\displacement - \vect{g}
    \right] &= \vect{0} && \text{on } (0, \tfinal) \times \partial\Omega_0.
    \label{eq:contact_boundary_omega_complementary_lagrangian}
  \end{empheq}
\end{subequations}
To solve this type of inequality boundary constraints, we employ a
\emph{primal-dual active set algorithm}~\cite{brunssen2007fast,
  hueber2005primal-dual}. This algorithm is introduced in
\cref{sec:numerics} and will be included in the numerical solution algorithm
interpreted as \emph{semismooth Newton algorithm}
\cite{hintermuller2002primal-dual}.

%%%%%%%%%%%%%%%%%%%%%%%%%%%%%%%%%%%%%%%%%%%%%%%%%%%%%%%%%%%%%%%%%%%%%%%%%%%%%%%%
%% End of 2-theory.tex
%%%%%%%%%%%%%%%%%%%%%%%%%%%%%%%%%%%%%%%%%%%%%%%%%%%%%%%%%%%%%%%%%%%%%%%%%%%%%%%%

  %% Numerics
  % !TeX encoding = UTF-8
% !TeX spellcheck = en_US

%%%%%%%%%%%%%%%%%%%%%%%%%%%%%%%%%%%%%%%%%%%%%%%%%%%%%%%%%%%%%%%%%%%%%%%%%%%%%%%%
%% 3-Numerics
%%%%%%%%%%%%%%%%%%%%%%%%%%%%%%%%%%%%%%%%%%%%%%%%%%%%%%%%%%%%%%%%%%%%%%%%%%%%%%%%

\section{Numerical Approach}
\label{sec:numerics}

This sections deals with the numerical treatment of the model equations.
Firstly, the normalization and mathematical problem is stated. Secondly, the
steps for solving the initial boundary value problem are stated including a
formulation for the space discretization with finite elements, time
discretization and the primal-dual active set algorithm as semismooth Newton
method. Finally, we incorporate the semismooth Newton method in the adaptive
space and time integration scheme \cite{castelli2021efficient} and propose the
numerical solution algorithm for the obstacle problem.

%%%%%%%%%%%%%%%%%%%%%%%%%%%%%%%%%%%%%%%%%%%%%%%%%%%%%%%%%%%%%%%%%%%%%%%%%%%%%%%%
\subsection{Problem Formulation}
\label{subsec:problem_statement}
%%%%%%%%%%%%%%%%%%%%%%%%%%%%%%%%%%%%%%%%%%%%%%%%%%%%%%%%%%%%%%%%%%%%%%%%%%%%%%%%

First, we improve the numerical stability by introducing a
nondimensionalization of the model equations. The \si{C}-rate specifies the
hours for the charging of the particle. We use the cycle time
$t_\text{cycle}= 1 / \text{C-rate}$ for the time scale, the particle
radius~$L_0$ in the Lagrangian frame for the spatial scale and the maximal
concentration~$c_{\max}$ as reference concentration. All dimensionless
variables are collected in \cref{tab:normalization}.
The dimensionless number $\tilde{\texttt{E}}_\text{H}$ is used to relate the
mechanical energy scale to the chemical energy scale, while the dimensionless
\emph{Fourier number}~$\Fo$ is used to relate the diffusion time scale to the
process time scale. From now, these dimensionless quantities are considered for
the model equations, however, we neglect accentuation to improve readability.

\begin{table}[t]
  \centering
  \caption{Dimensionless variables of the used model equations.}
  \label{tab:normalization}
  \scalebox{1.0}{
    \begin{tabular}[htbp]{@{}lllll@{}}
      \toprule
      $\tilde{t} = t / t_\text{cycle}$ & $\tilde{\vect{X}}_0 = \vect{X}_0 /
      L_0$ & $\tilde{\vect{u}} = \vect{u} / L_0 \qquad $ &
      $\tilde{c} =
      c /
      c_{\max} \qquad $ & $\tilde{v} = vc_{\max}$ \\
      $\tilde{U}_{\text{OCV}} = FU_{\text{OCV}} / R_\text{gas} T$ &
      $\tilde{\texttt{E}}_\text{H} = \EH / R_\text{gas} T c_{\max}$&
      \multicolumn{3}{c}{$\tilde{N}_{\text{ext}} =
        N_{\text{ext}}t_\text{cycle}/ L_0
        c_{\max} $ \qquad $\Fo = Dt_\text{cycle}/ L_0^2$} \\
      \bottomrule
    \end{tabular}
  }
\end{table}

For a general mathematical problem formulation, we follow the approach of
\cite{castelli2021efficient} and solve our set of equations for the
concentration $c$, the chemical potential $\mu$ and the displacement
$\vect{u}$. The mixed formulation of the Cahn--Hilliard-type equations would
allow to easily integrate the interfacial energy density for accounting phase
separation. The deformation gradient $\tens{F}$, the strain tensor
$\tens{E}_\el$ as well as the stress tensors $\tens{P}$ and
$\boldsymbol{\sigma}$ are calculated via the concentration $c$ and the
displacement~$\vect{u}$.

The dimensionless initial boundary value problem with inequality boundary
conditions is given as: let $\tfinal >0$ be the final simulation time and
$\Omega_0 \subset \IR^d$ a representative bounded electrode particle as
reference configuration with dimension~$d = 3$.
Find the normalized concentration~$c\colon[0, \tfinal] \times
\overline{\Omega}_0 \rightarrow [0,1]$, the chemical potential~$\mu\colon[0,
\tfinal] \times \overline{\Omega}_0 \rightarrow \IR$ and the
displacement~$\vect{u}\colon[0, \tfinal] \times \overline{\Omega}_0 \rightarrow
\IR^d$ satisfying

\begin{subequations}
  \label{eq:problem_boundary_constraints}
  \begin{empheq}[left=\empheqlbrace]{alignat=4}
    \partial_t c &= \minus \divgL\vect{N}(c, \grad_0\mu, \grad_0 \vect{u}) &&
    \qquad && \text{in } (0, \tfinal) \times
    \OmegaL,\label{eq:problem_boundary_constraintsa}\\
    \mu &= \partial_c \psi(c, \grad_0 \vect{u}) && && \text{in } (0, \tfinal)
    \times \OmegaL,\label{eq:problem_boundary_constraintsb}\\
    \vect{0} &= \minus \divgL \tens{P}(c, \grad_0\vect{u}) && && \text{in } (0,
    \tfinal) \times \OmegaL,\label{eq:problem_boundary_constraintsc}\\
    \vect{N} \cdot \normalL &= N_\text{ext} && && \text{on } (0, \tfinal)
    \times \partial \OmegaL,\label{eq:problem_boundary_constraintsd}\\
    \minus \tens{P} \cdot \normalL &\geq \vect{0} && && \text{on } (0, \tfinal)
    \times \partial\OmegaL,\label{eq:problem_boundary_constraintse}\\
    \displacement - \vect{g} &\leq \vect{0} &&  &&\text{on } (0, \tfinal)
    \times \partial\OmegaL,\label{eq:problem_boundary_constraintsf}\\
    \left[\minus \tens{P} \cdot \normalL\right] \left[\displacement-\vect{g}
    \right] &= \vect{0} && && \text{on } (0, \tfinal) \times
    \partial\OmegaL,\label{eq:problem_boundary_constraintsg}\\
    c(0, \something) &= c_0 && && \text{in } \OmegaL
    \label{eq:problem_boundary_constraintsh}
  \end{empheq}
\end{subequations}
and an initial condition $c_0$ that is consistent with the boundary conditions.
In case of lithiation we have a positive sign for the external lithium flux
$N_\text{ext}$ and in case of delithiation a negative sign. With appropriate
boundary conditions for the displacement, rigid body motions are excluded.
Note that the original formulation for the chemical deformation gradient
$\tens{F}_\ch$ is derived for the three-dimensional case, but all variables and
equations are mathematically valid in lower dimensions as well.

%%%%%%%%%%%%%%%%%%%%%%%%%%%%%%%%%%%%%%%%%%%%%%%%%%%%%%%%%%%%%%%%%%%%%%%%%%%%%%%%
\subsection{Numerical Solution Procedure}
\label{subsec:num_sol_approach}
%%%%%%%%%%%%%%%%%%%%%%%%%%%%%%%%%%%%%%%%%%%%%%%%%%%%%%%%%%%%%%%%%%%%%%%%%%%%%%%%

In this subsection, we present all details for the numerical solution of our
model equations: the space and time discretization of the initial boundary
value problem~\eqref{eq:problem_boundary_constraints}, the interpretation of the
primal-dual active set algorithm as semismooth Newton method and finally the
proposed adaptive solution algorithm.

%%%%%%%%%%%%%%%%%%%%%%%%%%%%%%%%%%%%%%%%%%%%%%%%%%%%%%%%%%%%%%%%%%%%%%%%%%%%%%%%
\subsubsection{Weak Formulation}
\label{subsubsec:weak_formulation}
%%%%%%%%%%%%%%%%%%%%%%%%%%%%%%%%%%%%%%%%%%%%%%%%%%%%%%%%%%%%%%%%%%%%%%%%%%%%%%%%

For the spatial discrete formulation, we derive the weak formulation of
\cref{eq:problem_boundary_constraints}. We define the function space
$\vect{V}^*\coloneqq H_*^1(\OmegaL; \IR^d)$ which includes appropriate boundary
constraints for the displacement considering possible boundary conditions
without contact. These displacement boundary constraints are stated for the
precise application case in \cref{sec:results}. Furthermore, we declare
the boundary on $\Omega_0$ as $\Gamma_\possibleDOF \coloneqq \ptl \Omega_0$ for
the potential contact zone to be in contact with the obstacle. For the
definitions of the scalar products, see \cref{app:tensor_analysis}.

The weak solution can be derived from a minimization problem on a convex set,
compare, e.g., \cite[Chapter~1.2]{kornhuber1997adaptive}
or~\cite{haslinger1980contact, boieri1987existence, kornhuber2001adaptive}, or
equivalently from a variational inequality, e.g.
\cite[Chapter~1.11]{friedman1982variational},
\cite[Chapter~2.1]{hlavacek1988solution},
\cite[Chapter~II.6]{kinderlehrer2000introduction},
\cite[Chapter~1.2]{kornhuber1997adaptive} or \cite{kornhuber2001adaptive}.
Multiplying with test functions, integration over the reference
domain~$\Omega_0$ and integration by parts, we state the weak formulation with a
variational inequality in the third equation: find the solutions $c, \mu,
\vect{u}$ with $c, \mu \in V \coloneqq H^1(\OmegaL)$, $\partial_t c \in
L^{2}(\OmegaL)$ and $\displacement \in \vect{V}^+ \coloneqq \left\{
\displacement \in \vect{V}^* \,:\, \displacement \leq \vect{g} \text{ on }
\Gamma_\possibleDOF \right\}$ such that
\begin{subequations}
  \label{eq:weak_chemical_mechanical_contact}
  \begin{empheq}[left=\empheqlbrace]{alignat=2}
    \label{eq:weak_chemical_mechanical_contacta}
    \left( \varphi,\partial_t c \right) &= \minus \Bigl( {m}(c,
    \grad_0\vect{u}) \gradL \varphi, \gradL\mu \Bigr) - \left(\varphi,
    {N}_\text{ext}\right)_{\Gamma_\possibleDOF}, \vspace{0.2cm}\\
    \label{eq:weak_chemical_mechanical_contactb}
    0 &= \minus \left(\varphi, \mu \right) + \bigl( \varphi, \partial_c
    \psi_\ch(c) + \partial_c \psi_\el(c, \grad_0\vect{u}) \bigr),
    \vspace{0.2cm}\\
    \label{eq:weak_chemical_mechanical_contactc}
    {0} &\leq \bigl( \gradL \left(\vect{\xi}-\displacement\right),
    \tens{P}(c, \grad_0 \vect{u})\bigr)
  \end{empheq}
\end{subequations}
for all test functions $\varphi \in V$ and $\vect{\xi} \in {\vect{V}}^+$. For
the formulation as saddle point problem, we follow
\cite{hueber2013contact, hueber2005primal-dual}, \texttt{deal.II} tutorial
step-41 in \cite{arndt2021deal-ii} and \cite{frohne2016efficient} and introduce
the Lagrange multiplier~$\vect{\lambda} \coloneqq \minus \boldsymbol{\sigma}
\cdot \normal = \minus \tens{P} \cdot \normalL$.

%%%%%%%%%%%%%%%%%%%%%%%%%%%%%%%%%%%%%%%%%%%%%%%%%%%%%%%%%%%%%%%%%%%%%%%%%%%%%%%%
\subsubsection{Space Discretization}
\label{subsubsec:space_discretization}
%%%%%%%%%%%%%%%%%%%%%%%%%%%%%%%%%%%%%%%%%%%%%%%%%%%%%%%%%%%%%%%%%%%%%%%%%%%%%%%%

For the spatial discretization we choose a computational domain $\Omega_h$ that
approximates the particle geometry $\Omega_0$ by a polytop. To approximate the
curved boundaries, we choose the isoparametric Lagrangian finite element
method~\cite[Chapter III \S{}2]{braess2007finite} on an admissible
mesh~${\mathcal{T}_n}$. We define the finite dimensional Lagrangian finite
element subspaces with the basis functions
\begin{subequations}
  \begin{alignat}{2}
    &{V}_h &&= \text{span} \{ \varphi_i \,:\, i=1,\dots, N \} \subset {V},\\
    &\vect{V}_h^+ &&= \text{span} \{ \vect{\xi}_j \,:\, j=1,\dots, d N \}
    \subset \vect{V}^+,\\
    &\vect{\Lambda}_h &&= \text{span}\{ \vect{\zeta}_k \,:\, k=1,\dots, d
    N_\Lambda
    \} \subset \vect{\Lambda},
  \end{alignat}
\end{subequations}
where $N$ denotes the number of degrees of freedom (DOFs) of the space ${V}_h$
and $N_\Lambda$ denotes the total number of nodes of the potential contact
zone $\Gamma_\possibleDOF$. For more details of the discretization, especially
the discretization of the Lagrange multiplier space $\vect{\Lambda}_h$, we
refer to
\cite{boieri1987existence, ben-belgacem1999extension, brunssen2007fast,
  hild2000numerical, wohlmuth2000mortar} and the references therein.

We now seek the discrete solutions for the concentration $\cD\colon [0,
\tfinal] \rightarrow \{ {V}_h: \cD \in [0,1] \}$, the chemical
potential~$\muD\colon [0, \tfinal] \rightarrow {V}_h$, the displacement
$\displacementD\colon[0, \tfinal]  \rightarrow \vect{V}_h^*$ and the Lagrange
multiplier~$\vect{\lambda}_h\colon[0, \tfinal] \rightarrow \vect{\Lambda}_h$ of
the spatial discrete saddle point problem
of~\cref{eq:weak_chemical_mechanical_contact}.

In a next step, we want to add the finite element ansatz. Therefore, we
represent the discrete solution variables with the basis functions given by
\begin{subequations}
  \begin{alignat}{2}
    &\cD(t, \vect{X}_0) = \sum_{i=1}^N c_i(t)\varphi_i(\vect{X}_0), \qquad
    &&\muD(t, \vect{X}_0) = \sum_{i=1}^N \mu_i(t)\varphi_i(\vect{X}_0),\\
    &\displacementD(t, \vect{X}_0) = \sum_{j=1}^{d N}
    u_j(t)\vect{\xi}_j(\vect{X}_0), \qquad
    &&\vect{\lambda}_h(t, \vect{X}_0) = \sum_{k=1}^{d N_{\Lambda}}
    \lambda_k(t)\vect{\zeta}_k(\vect{X}_0).
  \end{alignat}
\end{subequations}
For the vector valued finite dimensional subspace $\vect{V}_h^+ =
\text{span}\big\{ \vect{\xi}_j \,:\, j=1,\dots, d N \big\}$ and equivalent for
$\vect{\Lambda}_h$, we note~${\xi}_j$ as the scalar basis function, which is
the nonzero entry of the basis function vector~$\vect{{\xi}_j}$ of node $j$.
To simplify our notation, we use the same symbol for a function in
$\vect{V}_h^+$ and $\vect{\Lambda}_h$ as for its algebraic representation in
terms of the nodal basis. For the concentration and the chemical potential, we
use $\vect{c}_h$ and $\vect{\mu}_h$ for the algebraic representation.

Following \cite{brunssen2007fast, hueber2005primal-dual, hueber2005priori} the
biorthogonality of the basis functions has following property:
\begin{align}
  \label{eq:biortho}
  \int_{\Gamma_{\possibleDOF}} \vect{\xi}_j  \cdot \vect{\zeta}_k \de \vect{S}_0
  = \delta_{j,k} \int_{\Gamma_{\possibleDOF}} {\xi}_j \de \vect{S}_0
\end{align}
for all $j=1,\dots, d N$ and $k=1,\dots, d N_\Lambda$. For \cref{eq:biortho} we
suppose that the basis function $\vect{\xi}_j$ and the basis function
$\vect{\zeta}_j$ with the same index $j$ are associated to the same node on
$\Gamma_{\possibleDOF}$. The Kronecker symbol $\delta_{j,k}$ can be interpreted
as follows
\cite{brunssen2007fast}:
\begin{align}
  \delta_{j,k} = \left\{\begin{array}{ll}
    1, & \text{node $j$ coincides with potential contact node $k$},
    \vspace{0.2cm}\\
    0, & \text{otherwise.}
  \end{array}\right.
\end{align}

At the end of the spatial discretization process, we want to formulate our
problem as a discrete nonlinear differential algebraic equation~(DAE) before we
perform the time discretization in the next
\cref{subsubsec:time_discretization}. We therefore have a closer look on the
algebraic representation of our discrete weak formulation, in particular of the
momentum balance equation and of the displacement inequality at the boundary
condition \eqref{eq:problem_boundary_constraintsf}.

Let $\vect{u}_h$ and $\vect{\lambda}_h$ be the solution of the discrete
variational inequality. Then, we have the algebraic representation of the
discrete weak formulation of \cref{eq:problem_boundary_constraintsc} as
\begin{align}
  \label{eq:discrete_momentum_balance}
  {\tens{P}}_h(c_h, \grad_0{\vect{u}}_h) + \tens{B}_h {\vect{\lambda}}_h &=
  \vect{0}
\end{align}
with the nonlinear vector
\begin{align}
  \label{eq:matrix_P_h}
  \tens{P}_h(c_h, \grad_0{\vect{u}}_h) = \left[\int_{\OmegaL}
  \gradL\vect{\xi}_j \!:\! \tens{P}\bigl(\cD, \gradL\vect{u}_h\bigr)
  \de \vect{X}_0\right]_{j = 1, \dots, dN}
\end{align}
and the matrix
\begin{align}
  \label{eq:matrix_B_h}
  \tens{B}_h = \left[\int_{\Gamma_{\possibleDOF}} \vect{\xi}_j \cdot
  \vect{\zeta}_k\de \vect{S}_0\right]_{j=1,\dots, d N,\ k=1,\dots, d N_\Lambda}.
\end{align}
With an appropriate node numbering, $\tens{B}_h$ can be written as
$\tens{B}_h = \left( \vect{0}, \tens{D}_h \right)^\trp$. Due to the
biorthogonality in \cref{eq:biortho}, the diagonal matrix $\tens{D}_h$ has the
entries
\begin{align}\label{eq:diagonal_help}
  \left({D}_h\right)_{j,k} = \int_{\Gamma_{\possibleDOF}} \vect{\xi}_j \cdot
  \vect{\zeta}_k \de \vect{S}_0 = \delta_{j,k}  \int_{\Gamma_{\possibleDOF}}
  {\xi}_k \de \vect{S}_0
\end{align}
for all $j,k=1,\dots, d N_{\Lambda}$. Now, we define two sets of indices: all
degrees of freedom on $\Gamma_{\possibleDOF}$ with $\possibleDOF$ (all
potential contact nodes) and all other nodes with $\otherDOF$.

Consider now the weaker integral condition for the strong pointwise
non-penetration condition of \cref{eq:problem_boundary_constraints} for the
discrete contact conditions for all $p \in \possibleDOF$:
\begin{align}
  \label{eq:weak_non_penetration}
  \int_{\Gamma_{\possibleDOF}} \displacement_{h} \cdot \vect{\zeta}_p
  \de\vect{S}_0 \leq \int_{\Gamma_{\possibleDOF}} \vect{g}_h \cdot
  \vect{\zeta}_p\de \vect{S}_0 \eqqcolon \hat{{g}}_{p} \quad
  \Longleftrightarrow \quad
  \int_{\Gamma_{\possibleDOF}} \displacementt_{p} {\xi}_p{\zeta}_p \de\vect{S}_0
  \leq \hat{{g}}_{p},
\end{align}
where $\displacementt_p$ is the scalar coefficient of the discrete vector
${\vect{u}}_h$ of DOF~$p$ and ${\vect{g}}_h$ is an appropriate approximation of
$\vect{g}$ on $\Gamma_{\possibleDOF}$. Next, we rewrite
\cref{eq:weak_non_penetration} with the help of \cref{eq:diagonal_help} for the
algebraic representation of the weak non-penetration condition
$\hat{\displacementt}_p \coloneqq \left({D}_h\right)_{p,p} \displacementt_p
\leq \hat{{g}}_{p}$ for all $p \in \possibleDOF$. Then, we can rewrite the
condition for the Lagrange multiplier in the same way to get
$\hat{{\lambda}}_{p} \coloneqq \left({D}_h\right)_{p,p} {\lambda}_p$
with the same definition for ${\lambda}_p$ as for ${\displacementt_p}$ for all
nodes $p \in \possibleDOF$.

Finally, the discrete algebraic form of the contact problem of
\cref{eq:problem_boundary_constraints} is given by:
\begin{subequations}
  \label{eq:discrete_contact_KKT}
  \begin{alignat}{1}
    &\qquad \quad {\tens{P}}_h(c_h, \grad_0{\vect{u}}_h) + \tens{B}_h
    \vect{\lambda}_h = \vect{0},\label{eq:discrete_contact}\\
    &\hat{{\displacementt}}_p \leq \hat{{g}}_p, \qquad
    \hat{{\lambda}}_{p} \geq
    {0}, \qquad \hat{{\lambda}}_{p}(\hat{{\displacementt}}_p -
    \hat{{g}}_{p}) =
    {0} \label{eq:discrete_KKT}
  \end{alignat}
\end{subequations}
for all nodes $p\in \possibleDOF$. \cref{eq:discrete_KKT} can also be
identified as discrete KKT conditions of a constrained optimization problem for
inequality constraints \cite{hueber2005primal-dual}.

Next, we use a reformulation of the three equations of \cref{eq:discrete_KKT}
based on the nonlinear complementarity problem (NCP) function
\begin{align}
  \mathcalC(a, b) := b -\max\left( b + \alpha a, 0 \right),
  \qquad \forall a,b \in \IR,
\end{align}
and $\alpha > 0$ arbitrarily fixed \cite{hintermuller2003semismooth}. The
following equivalence is true \cite{hintermuller2002primal-dual,
  hintermuller2003semismooth}:
\begin{align}
  \mathcalC(a,b) = 0 \quad \Longleftrightarrow \quad a \leq 0, \:\: b \geq
  0,\:\: ab = 0.
\end{align}
Applied to \cref{eq:discrete_KKT} it follows:
\begin{align}
  \label{eq:definition_of_mathcalC}
  {\mathcalC}(\hat{\displacementt}_p, \hat{{\lambda}}_{p}) =
  \hat{{\lambda}}_{p} -{\max}\left( \hat{{\lambda}}_{p} + \alpha
  \left(\hat{{\displacementt}}_p-\hat{{g}}_{p}\right)\!, {0} \right) = {0}
\end{align}
for all $p \in \possibleDOF$ and $\alpha > 0$. In total we rewrite
\cref{eq:discrete_contact_KKT} to
\begin{subequations}
  \label{eq:discrete_final}
  \begin{alignat}{2}
    &{\tens{P}}_h(c_h, \grad_0{\vect{u}}_h) + \tens{B}_h {\vect{\lambda}}_h
    &&= \vect{0},\label{eq:discrete_contact_final}\\
    &\boldsymbol{\mathcalC}_\possibleDOF({\vect{u}}_h, {\vect{\lambda}}_{h})
    &&= \vect{0} \label{eq:discrete_KKT_final}
  \end{alignat}
\end{subequations}
with the same definition for $\boldsymbol{\mathcalC}_\possibleDOF(\something,
\something)$ in each component as in \cref{eq:definition_of_mathcalC} for all
$p \in
\possibleDOF$, compare also~\cref{app:tensor_analysis}.

\begin{figure}[t]
  \centering
  \includegraphics[width = 0.38\textwidth,
  page=1]{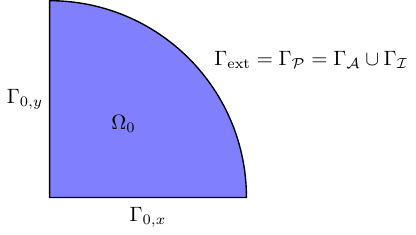}
  \caption{All boundary parts for the example of a two-dimensional quarter
    disk domain $\Omega_0$ split up into two artificial
    boundaries
    $\Gamma_{0,y}$ and $\Gamma_{0,x}$ with additional Dirichlet constraints for
    the displacement~$\vect{u}$ and the potential contact
    boundary~$\Gamma_\possibleDOF$
    subdivided into the active contact boundary~$\Gamma_\activeset$ and the
    inactive boundary~$\Gamma_\inactiveset$.}
  \label{fig:2D_quarterdisk_different_sets}
\end{figure}

Considering for example an quarter shaped obstacle like in
\cref{fig:2D_quarterdisk_different_sets}, the physical boundary of~$\OmegaL$
can be defined as~$\Gamma_\text{ext}$. $\Gamma_{0,y}$ and $\Gamma_{0,x}$ are two
artificial boundaries with appropriate Dirichlet boundary conditions for the
displacement~$\vect{u}$. $\Gamma_\text{ext}$ is the potential contact
boundary~$\Gamma_\possibleDOF$ and splits into two parts, the active contact
boundary~$\Gamma_\activeset$ and the inactive boundary~$\Gamma_\inactiveset$.

Collecting all time-dependent solution variables in a vector-valued function
\begin{align}
  \label{eq:discrete_solution_vector}
  {\vect{y}}\colon[0, \tfinal] \rightarrow \IR^{(2+d)N+dN_\Lambda}, \quad t
  \mapsto \vect{y}(t) =
  \begin{pmatrix}
    {\vect{c}}_h \\
    {\vect{\mu}}_h \\
    {\vect{u}}_h \\
    {\vect{\lambda}}_h,
  \end{pmatrix}
\end{align}
we can state our spatial discrete problem of the saddle point formulation
of~\cref{eq:weak_chemical_mechanical_contact} as general nonlinear DAE: Find
$\vect{y}\colon [0, \tfinal] \rightarrow \IR^{(2+d)N+dN_\Lambda}$ satisfying
\begin{align}
  \label{eq:resulsting_dae}
  \tens{M} \ptl_t \vect{y} = \vect{f}(t, \vect{y}) \qquad \text{for } t\in (0,
  \tfinal], \qquad \vect{y}(0) = \vect{y}^0.
\end{align}
On the left side, the system mass matrix $\tens{M}$ has only one nonzero-block
entry $\tens{M}_{h} = \left[(\varphi_i, \varphi_j)\right]_{ij}$ and is
therefore singular. $\tens{M}_{h}$ identifies the mass matrix of the finite
element space $V_h$. On the right side, $\vect{f}\colon [0, \tfinal] \times
\IR^{(2+d)N+dN_\Lambda} \rightarrow \IR^{(2+d)N+dN_\Lambda}$ is given according
to the algebraic formulation in \cref{eq:discrete_final} by
\begin{align}
  \label{eq:right_hand_side_dae}
  \vect{f}(t, \vect{y}) \coloneqq
  \begin{pmatrix}
    \minus \tens{K}_{m}(c_h, \grad_0 \vect{u}_h){\vect{\mu}}_h -
    \vect{N}_\text{ext} \\
    \minus \tens{M}_h {\vect{\mu}}_h + \vect{\Psi}_\ch(c_h) +
    \vect{\Psi}_\el(c_h, \grad_0 \vect{u}_h) \\
    {\tens{P}}_h(c_h, \grad_0{\vect{u}}_h) + \tens{B}_h {\vect{\lambda}}_h \\
    \boldsymbol{\mathcalC}_\possibleDOF({\vect{u}}_h,{\vect{\lambda}}_h)
  \end{pmatrix}
\end{align}
with the relation of $\vect{y}$ to the solution variables as in
\cref{eq:discrete_solution_vector}. The quantities in the definition of
$\vect{f}$ are given by the mass matrix $\tens{M}_h$,
the stiffness matrix $\tens{K}_{m}(c_h, \grad_0 \vect{u}_h) =
\left[\left({m}(c_h, \grad_0 \vect{u}_h)\grad_0 \varphi_i,
\grad_0\varphi_j\right)\right]_{i,j}$, the vectors for the nonlinearities
$\vect{\Psi}_\ch(c_h) = \left[\left(\varphi_i, \ptl_c \psi_\ch(c_h)
\right)\right]_{i}$ and $\vect{\Psi}_\el(c_h, \grad_0 \vect{u}_h) =
\left[\left(\varphi_i, \ptl_c \psi_\el(c_h, \grad_0 \vect{u}_h)
\right)\right]_{i}$, as well as the boundary
condition~\mbox{$\vect{N}_\text{ext} =
  \left[\left(\varphi_i, N_\text{ext}\right)_{\Gamma_\text{ext}}\right]_i$.} For
the block
mass matrix~$\tens{B}_h$, the nonlinear vector~${\tens{P}}_h(c_h,
\grad_0{\vect{u}}_h)$ and the
vector~$\boldsymbol{\mathcalC}_\possibleDOF({\vect{u}}_h,{\vect{\lambda}}_h)$,
we refer to the definitions
in~\cref{eq:matrix_B_h,eq:matrix_P_h,eq:discrete_KKT_final}, respectively.

%%%%%%%%%%%%%%%%%%%%%%%%%%%%%%%%%%%%%%%%%%%%%%%%%%%%%%%%%%%%%%%%%%%%%%%%%%%%%%%%
\subsubsection{Time Discretization}
\label{subsubsec:time_discretization}
%%%%%%%%%%%%%%%%%%%%%%%%%%%%%%%%%%%%%%%%%%%%%%%%%%%%%%%%%%%%%%%%%%%%%%%%%%%%%%%%

For the temporal discretization, we follow the approach
in~\cite{castelli2021efficient} by using a variable-step, variable-order
algorithm~\cite{reichelt1997matlab, shampine1997matlab, shampine1999solving,
  shampine2003solving}. This approach seems reasonable since the
DAE~\eqref{eq:resulsting_dae} can be treated in a similar way to a stiff
ordinary differential equation. The algorithm adaptively changes the time step
size $\tau_n > 0$ and the order by an error control.

This leads to the space and time discrete problem: find the discrete
solution $\vect{y}^{n+1} \approx \vect{y}(t_{n+1})$ satisfying
\begin{align}
  \label{eq:space_and_time_discretization}
  \alpha_{k_n} \tens{M}\left(\vect{y}^{n+1} - \vect{\Phi}^n\right) = \tau_n
  \vect{f}\left(t_{n+1}, \vect{y}^{n+1}\right)
\end{align}
to advance one time step from $t_n$ to $t_{n+1} = t_n + \tau_n$.
$\vect{\Phi}^n$ is defined by the solutions on the former time steps
$\vect{y}^n, \dots, \vect{y}^{n-k}$ and a constant $\alpha_{k_n} > 0$ depending
on the selected order $k_n$ at time $t_n$
\cite[Section~2.3]{shampine1997matlab}. Because of the time-dependent Neumann
boundary condition $N_\text{ext}$, the vector $\vect{f}$ does also explicitly
depend on the time~$t$.

%%%%%%%%%%%%%%%%%%%%%%%%%%%%%%%%%%%%%%%%%%%%%%%%%%%%%%%%%%%%%%%%%%%%%%%%%%%%%%%%
\subsubsection{The Primal-Dual Active Set Algorithm as Semismooth Newton Method}
\label{subsubsec:primal_dual_active_set_algorithm_as_semismooth_newton_method}
%%%%%%%%%%%%%%%%%%%%%%%%%%%%%%%%%%%%%%%%%%%%%%%%%%%%%%%%%%%%%%%%%%%%%%%%%%%%%%%%

The next step is to use an appropriate iterative solution scheme to handle the
NCP function in \cref{eq:discrete_KKT_final}. The primal-dual active set
algorithm is the strategy of choice, since it is an iterative approach to deal
with the condition in \cref{eq:discrete_KKT_final} and to predict the next
active and inactive set $\activeset_{k+1}^{n+1}$ and $\inactiveset_{k+1}^{n+1}$
for each new time step~$t_{n+1}$ \cite{brunssen2007fast, hueber2005primal-dual,
  hintermuller2002primal-dual}.

For the moment, we consider only \cref{eq:discrete_contact_final} and
\cref{eq:discrete_KKT_final}, which are the relevant parts for the contact
inequality condition. To compute the new active and inactive set
$\activeset_{k+1}^{n+1}$ and $\inactiveset_{k+1}^{n+1}$ on the potential
boundary $\Gamma_{\possibleDOF}$ in a new time
step~$t_{n+1}$ we will use the \emph{primal-dual active set algorithm},
compare~\cite{brunssen2007fast, hueber2005primal-dual,
  hintermuller2002primal-dual}.

Since we have to linearize \cref{eq:space_and_time_discretization} anyway to
handle the nonlinear algebraic system via the Newton--Raphson method, we can
also use the interpretation of the primal-dual active set algorithm as
semismooth Newton method, compare \cite{frohne2016efficient,
  hintermuller2002primal-dual, hintermuller2003semismooth} and \texttt{deal.II}
tutorials step-41 and step-42 in \cite{arndt2021deal-ii}. The locally
superlinear convergence and global convergence results are shown
in~\cite{hintermuller2002primal-dual}. Since $\mathcalC(\something,
\something)$ in \cref{eq:discrete_KKT_final} is not differentiable, Newton
techniques for solving \cref{eq:discrete_final} have to be applicable with
generalizations of the derivative of a function. These methods are named
\textit{generalized Newton methods}, see e.g. \cite{hintermuller2003semismooth}
and the references therein. In this paper we propose the \textit{semismooth
  Newton method} following \cite{hintermuller2003semismooth,
  hintermuller2002primal-dual}. The semismoothness characteristics of the
$\max$-operator in \cref{eq:discrete_KKT_final} lead to local convergence
properties of the semismooth Newton method.

So we need a linearization of a function $\mathcalC(\something, \something)$,
being not classically differentiable. As a replacement, we can use the concept
of slant differentiability, compare~\cite{hintermuller2003semismooth}.
$\mathcalC(\something, \something)$ is slantly differentiable with
\begin{subequations}
  \begin{alignat}{3}
    &\FP{}{\displacementt_p} \mathcalC(\hat{\displacementt}_{p},
    \hat{{\lambda}}_{p})
    &&= \left\{\begin{array}{ll}
      \minus \alpha \left({D}_h\right)_{p,p},
      & \text{ if } \hat{{\lambda}}_{p} +
      \alpha (\hat{{\displacementt}}_{p}-\hat{{g}}_{p}) > 0, \vspace{0.2cm}\\
      0,
      & \text{ if } \hat{\lambda}_{p} + \alpha
      (\hat{\displacementt}_{p}-\hat{{g}}_{p}) \leq 0,
    \end{array}\right.\\
    &\FP{}{{\lambda}_p} \mathcalC(\hat{\displacementt}_{p},
    \hat{{\lambda}}_{p})
    &&= \left\{\begin{array}{ll} 0,
      & \hspace*{1.2cm} \text{ if } \hat{\lambda}_{p} + \alpha
      (\hat{\displacementt}_{p}-\hat{{g}}_{p}) > 0, \vspace{0.2cm}\\
      1,
      & \hspace*{1.2cm} \text{ if } \hat{\lambda}_{p} + \alpha
      (\hat{\displacementt}_{p}-\hat{{g}}_{p}) \leq 0.
    \end{array}\right.
  \end{alignat}
\end{subequations}
Recall the definition of the potential contact nodes $\possibleDOF = \activeset
\cup \inactiveset$, the active and inactive sets, respectively, and the
remaining nodes $\otherDOF$. The set of all nodes is defined via $\alldofs =
\possibleDOF \cup \otherDOF$. Moreover, we define for the Jacobian matrix the
derivatives for the nonlinear part:
\begin{align}
  {}_{z}\tens{A}_h &= \left[\left(\grad_0 \vect{\xi}_k, \ptl_z \tens{P}(c_h^k,
  \grad_0 \vect{u}_h^k) \varphi_{i} \right) \right]_{k,i}, \qquad
  {}_{\tens{G}}\tens{A}_h = \left[\left(\grad_0 \vect{\xi}_k, \ptl_{\tens{G}}
  \tens{P}(c_h^k, \grad_0 \vect{\xi}_l) \right) \right]_{k,l},
\end{align}
for $i= 1,\dots, N$ and $k,l = 1,\dots, dN$, with the derivative regarding the
first scalar valued quantity and the second tensor valued quantity,
respectively. With the partitions of the matrices~${}_{z}\tens{A}_h$,
${}_{\tens{G}}\tens{A}_h$ and $\tens{B}_h$ as well as the vectors $\vect{u}_h$
and $\vect{\lambda}_h$ for the different sets, a semismooth Newton step for
\cref{eq:discrete_final} has the form
\begin{align}\label{eq:semismooth_newton_update}
  %&
  \begin{pmatrix}
    {}_{z}\tens{A}_{\otherDOF\alldofs} &
    {}_{\tens{G}}\tens{A}_{\otherDOF\otherDOF} &
    {}_{\tens{G}}\tens{A}_{\otherDOF\inactiveset_k} &
    {}_{\tens{G}}\tens{A}_{\otherDOF\activeset_k} & \tens{0} & \tens{0}\\
    {}_{z}\tens{A}_{\inactiveset_k\alldofs} &
    {}_{\tens{G}}\tens{A}_{\inactiveset_k\otherDOF} &
    {}_{\tens{G}}\tens{A}_{\inactiveset_k\inactiveset_k} &
    {}_{\tens{G}}\tens{A}_{\inactiveset_k\activeset_k} &
    \tens{D}_{\inactiveset_k} & \tens{0}\\
    {}_{z}\tens{A}_{\activeset_k\alldofs} &
    {}_{\tens{G}}\tens{A}_{\activeset_k\otherDOF} &
    {}_{\tens{G}}\tens{A}_{\activeset_k\inactiveset_k} &
    {}_{\tens{G}}\tens{A}_{\activeset_k\activeset_k} & \tens{0} &
    \tens{D}_{\activeset_k}\\
    \tens{0} & \tens{0} & \tens{0} & \tens{0} & \tens{Id}_{\inactiveset_k} &
    \tens{0}
    \\
    \tens{0} & \tens{0} & \tens{0} & \minus \alpha \tens{D}_{\activeset_k} &
    \tens{0} & \tens{0}
  \end{pmatrix}
  \begin{pmatrix}
    \delta \vect{c}_{\alldofs}\\
    \delta\displacement_{\otherDOF}\\
    \delta\displacement_{\inactiveset_k}\\
    \delta\displacement_{\activeset_k}\\
    \delta\vect{\lambda}_{\inactiveset_k}\\
    \delta\vect{\lambda}_{\activeset_k}\\
  \end{pmatrix}
  = \minus \begin{pmatrix}
    {\tens{P}}_\otherDOF(c_h^k, \grad_0{\vect{u}}_h^k) \\
    {\tens{P}}_{\inactiveset_k}(c_h^k, \grad_0{\vect{u}}_h^k) +
    \left(\tens{B}_h \vect{\lambda}_h^k \right)_{\inactiveset_k}\\
    {\tens{P}}_{\activeset_k}(c_h^k, \grad_0{\vect{u}}_h^k) + \left(\tens{B}_h
    \vect{\lambda}_h^k \right)_{\activeset_k}\\
    \vect{\lambda}_{\inactiveset_k}^k\\
    \minus \alpha
    \left(\tens{D}_{\activeset_k}\displacement_{\activeset_k}^k
    -\vect{g}_{\activeset_k}\right)
    %\nonumber
  \end{pmatrix}
\end{align}
with the definition of \cref{eq:weak_non_penetration} for all active points
$p\in \activeset_k$, for $\vect{g}_{\activeset_k}$ and
$\left(\tens{B}_h\right)_\otherDOF = \tens{0}$. $\tens{Id}_{\inactiveset_k}$ is
the identity matrix of dimension $\text{card}\left(\inactiveset_k\right)$.

As a next step we have a closer look at different subsystems
of~\cref{eq:semismooth_newton_update}. The fourth row provides
\begin{align}
  \label{eq:recover_multiplier_inactive}
  \vect{\lambda}_{\inactiveset_k}^{k+1} = \vect{\lambda}_{\inactiveset_k}^{k} +
  \delta\vect{\lambda}_{\inactiveset_k} = \vect{\lambda}_{\inactiveset_k}^{k} -
  \vect{\lambda}_{\inactiveset_k}^{k} = \vect{0}
\end{align}
and the last one implies
\begin{align}\label{eq:displacement_new}
  \displacement_{\activeset_k}^{k+1} = \displacement_{\activeset_k}^{k} +
  \delta\displacement_{\activeset_k} =  \displacement_{\activeset_k}^{k}
  +\tens{D}_{\activeset_k}^{\minus 1}\vect{g}_{\activeset_k} -
  \displacement_{\activeset_k}^{k} = \tens{D}_{\activeset_k}^{\minus
    1}\vect{g}_{\activeset_k}.
\end{align}
\cref{eq:recover_multiplier_inactive} and \cref{eq:displacement_new} are
exactly the conditions of the active and inactive sets in the primal-dual active
set algorithm. Considering now the subsystem
of~\cref{eq:semismooth_newton_update} for the active set of the Lagrange
multiplier, we have
\begin{alignat}{2}
  \label{eq:recover_multiplier_active}
  &\vect{\lambda}_{\activeset_{k}}^{k+1} &&=
  \vect{\lambda}_{\activeset_{k}}^{k} + \delta\vect{\lambda}_{\activeset_{k}} =
  \minus \tens{D}_{\activeset_k}^{\minus 1} {\tens{P}}_{\activeset_k}(c_h^k,
  \grad_0{\vect{u}}_h^k) -\tens{D}_{\activeset_k}^{\minus
    1}\left({}_{\tens{G}}\tens{A}_{h} \delta\displacement
  \right)_{\activeset_k}
  - \tens{D}_{\activeset_k}^{\minus 1}\left({}_{z}\tens{A}_{h}\delta
  \vect{c}\right)_{\activeset_k}.
\end{alignat}
This means that the Lagrange multiplier only has to be computed on the active
set with the solutions~$\delta \vect{c}$ and $\delta \vect{u}$.

Let us now consider for a moment the two sets $\otherDOF$ and $\inactiveset$
together as $\otherandinactive$ since the subsystems
of~\cref{eq:semismooth_newton_update} for both sets are equal. The combination
of the two sets results in
\begin{align}\label{eq:combination_of_two_sets_together}
  \begin{pmatrix}
    {}_{z}\tens{A}_{\otherandinactive,\alldofs} &
    {}_{\tens{G}}\tens{A}_{\otherandinactive} \\
  \end{pmatrix}
  \begin{pmatrix}
    \delta \vect{c}_{\alldofs}\\
    \delta\displacement_{\otherandinactive}\\
  \end{pmatrix}
  = \minus \begin{pmatrix}
    {\tens{P}}_{\otherandinactive}(c_h^k, \grad_0{\vect{u}}_h^k) +
    {}_{\tens{G}}\tens{A}_{\otherandinactive,\activeset_k}
    \delta\displacement_{\activeset_{k}}
  \end{pmatrix}
\end{align}
to compute $\displacement_\otherandinactive^{k+1}$. In the end we solve the
reduced system for the Newton update
\begin{align}\label{eq:semismooth_newton_update_small}
  \begin{pmatrix}
    {}_{z}\tens{A}_{\otherDOF\alldofs} &
    {}_{\tens{G}}\tens{A}_{\otherDOF\otherDOF} &
    {}_{\tens{G}}\tens{A}_{\otherDOF\inactiveset_k} &
    {}_{\tens{G}}\tens{A}_{\otherDOF\activeset_k} \\
    {}_{z}\tens{A}_{\inactiveset_k\alldofs} &
    {}_{\tens{G}}\tens{A}_{\inactiveset_k\otherDOF} &
    {}_{\tens{G}}\tens{A}_{\inactiveset_k\inactiveset_k} &
    {}_{\tens{G}}\tens{A}_{\inactiveset_k\activeset_k}\\
    {}_{z}\tens{A}_{\activeset_k\alldofs} &
    {}_{\tens{G}}\tens{A}_{\activeset_k\otherDOF} &
    {}_{\tens{G}}\tens{A}_{\activeset_k\inactiveset_k} &
    {}_{\tens{G}}\tens{A}_{\activeset_k\activeset_k}
  \end{pmatrix}
  \begin{pmatrix}
    \delta \vect{c}_{\alldofs}\\
    \delta\displacement_{\otherDOF}\\
    \delta\displacement_{\inactiveset_k}\\
    \delta\displacement_{\activeset_k}
  \end{pmatrix}
  = \minus \begin{pmatrix}
    {\tens{P}}_{\otherDOF}(c_h^k, \grad_0{\vect{u}}_h^k)\\
    {\tens{P}}_{\inactiveset_k}(c_h^k, \grad_0{\vect{u}}_h^k)\\
    {\tens{P}}_{\activeset_k}(c_h^k, \grad_0{\vect{u}}_h^k)
  \end{pmatrix},
\end{align}
restricting the Newton update to zero for the degrees of freedom in the active
set and providing the correct boundary values to the new solution as
inhomogeneous Dirichlet boundary values. This can be done via $\vect{u}_h^k
\coloneqq P_{\activeset_{k+1}}(\displacement_{h}^{k})$ with the projection
\begin{align}
  \label{eq:displacement_projection}
  P_{\activeset_{k+1}}(\displacement_{h})_p \coloneqq
  \left\{\begin{array}{rl}
    \hat{\displacementt}_{p}, & \text{ if } p \not\in \activeset_{k+1},\\
    \hat{g}_{p}, & \text{ if } p \in \activeset_{k+1}.
  \end{array}\right.
\end{align}

After solving the total Newton system and computing the new solutions, the
Lagrange multiplier $\vect{\lambda}^{k+1}$ can be recovered via
\cref{eq:recover_multiplier_active} and \cref{eq:recover_multiplier_inactive}.

Following \cite{brunssen2007fast} we use
\begin{alignat}{2}
  \label{eq:recover_multiplier_active_used}
  &\vect{\lambda}_{\activeset_{k}}^{k+1} &&=
  \vect{\lambda}_{\activeset_{k}}^{k} + \delta\vect{\lambda}_{\activeset_{k}} =
  \minus \tens{D}_{\activeset_k}^{\minus 1}
  {\tens{P}}_{\activeset_k}(c_h^{k+1}, \grad_0{\vect{u}}_h^{k+1})
\end{alignat}
which has the same linearization as \cref{eq:recover_multiplier_active}. This
corresponds to an inexact strategy, compare Algorithm~3
in~\cite{brunssen2007fast}. Further, it reduces the computational effort of a
second nested loop compared to Algorithm~2 in~\cite{brunssen2007fast}. However,
the inexact case is an additional simplification of the applied algorithm used
in this paper. This further means that it has not been clarified whether the
superlinear convergence is retained \cite{brunssen2007fast}.

Finally, we have to formulate a semismooth Newton algorithm in one time step.
With the concept of the semismooth Newton algorithm we can also update our
DAE~\eqref{eq:resulsting_dae} and we can remove all parts related with the
Lagrange multiplier $\vect{\lambda}$. This results in a new definition of
$\tilde{\vect{y}}$, $\tilde{\vect{f}} \in \IR^{(2+d)N}$.
In the following, we consider only the updated system and omit again the
accentuation~$\tilde{\square}$. So we have to linearize the updated version of
the DAE~\eqref{eq:space_and_time_discretization}
to compute the Newton update.

%%%%%%%%%%%%%%%%%%%%%%%%%%%%%%%%%%%%%%%%%%%%%%%%%%%%%%%%%%%%%%%%%%%%%%%%%%%%%%%%
\subsubsection{Adaptive Solution Algorithm}
\label{subsubsec:adaptive_algorithm}
%%%%%%%%%%%%%%%%%%%%%%%%%%%%%%%%%%%%%%%%%%%%%%%%%%%%%%%%%%%%%%%%%%%%%%%%%%%%%%%%

After the linearization of the updated version of the
DAE~\eqref{eq:space_and_time_discretization} with the semismooth Newton method,
the Newton update is computed with a direct LU-decomposition. Keep in mind that
the number of iteration steps during the Newton method can be reduced with an
appropriate initialization. The starting values for the first time step are
given in \cref{subsec:simulation_setup} whereas during time integration a
predictor scheme is applied \cite{shampine1997matlab}.

For the space and time adaptive solution algorithm we follow Algorithm~1 in
\cite{castelli2021efficient}. Here, a temporal error
estimator~\cite{reichelt1997matlab, shampine1997matlab, shampine1999solving,
  shampine2003solving} and a spacial error estimator are applied. A gradient
recovery estimator is used for the spatial
regularity~\cite[Chapter~4]{ainsworth2000posteriori}. To mark the cells for
local coarsening and refinement, the parameters $\theta_{\mathrm{c}}$ and
$\theta_{\mathrm{r}}$ are used with a maximal
strategy~\cite{banas2008adaptive}. Finally, a mixed error control is applied
using the parameters $\text{RelTol}_t$, $\text{AbsTol}_t$, $\text{RelTol}_x$
and $\text{AbsTol}_x$. Further details can be found
in~\cite{castelli2021efficient}.

Combining the semismooth Newton method and the space and time adaptive
algorithm by~\cite[Algorithm~1]{castelli2021efficient}, we propose the
following concept:
\begin{center}
  \scalebox{0.85}
  {
    \begin{minipage}{1.05\linewidth}
      \setcounter{algorithm}{1}
      \begin{algorithm}[H]
        \caption{Semismooth Newton Method for Adaptive Obstacle Space Time
          Algorithm}
        \label{alg:cap}
        \begin{algorithmic}[1]
          \State
          Initialize $\activeset_1^{0}$ and $\inactiveset_1^{0}$ such that
          $\possibleDOF^0 = \activeset_1^0 \cup \inactiveset_1^0$ and
          $\activeset_1^0 \cap \inactiveset_1^0 = \emptyset$
          \While{$t_n < \tfinal$}
          \State
          Given ${\mathcal{T}_n}, \tau_n, k_n$ and $\sol^n, \dots,
          \sol^{n-k_n}$, set
          $k=1$
          \State
          Extrapolate $\sol^n, \dots, \sol^{n-k_n}$ to compute predictor
          $\sol^{(0),n+1}$ in $t_{n+1}$
          \While{not converged}
          \State
          Solve for Newton update $\delta\sol^{k,n}$ (introducing the contact
          condition as additional Dirichlet \newline
          \hspace*{1.35cm} boundary condition and set the Newton update to zero
          for DOFs in the active set)
          \State
          Compute $\sol^{k+1,n+1}$ = $\sol^{k,n} + \delta\sol^{k,n}$
          \State
          Recover $\vect{\lambda}^{k+1,n+1}$ and compute new
          $\activeset_{k+1}^{n+1}$ and $\inactiveset_{k+1}^{n+1}$
          \State
          Project $\sol^{k,n+1}$ according to $\activeset_{k+1}^{n+1}$ and
          update constraints
          \If{$\activeset_{k+1}^{n+1} = \activesetk^{n+1}$ and the Newton
            update norm is reduced appropriately}
          \State
          Exit inner while loop
          \ElsIf{Newton update norm is not reduced appropriately or maximal
            Newton iteration
            number is reached}
          \State
          Reduce time step size and go to Line 3
          \Else
          \State
          Update $k+1 \rightarrow k$ and go to Line 7
          \EndIf
          \EndWhile
          \State
          Advance time step via space and time algorithm
          \cite[Algorithm~1]{castelli2021efficient}
          \EndWhile
        \end{algorithmic}
      \end{algorithm}
    \end{minipage}
  }
\end{center}

%%%%%%%%%%%%%%%%%%%%%%%%%%%%%%%%%%%%%%%%%%%%%%%%%%%%%%%%%%%%%%%%%%%%%%%%%%%%%%%%%
%%% End of 3-numerics.tex
%%%%%%%%%%%%%%%%%%%%%%%%%%%%%%%%%%%%%%%%%%%%%%%%%%%%%%%%%%%%%%%%%%%%%%%%%%%%%%%%%

  %% Results
  % !TeX encoding = UTF-8
% !TeX spellcheck = en_US

%%%%%%%%%%%%%%%%%%%%%%%%%%%%%%%%%%%%%%%%%%%%%%%%%%%%%%%%%%%%%%%%%%%%%%%%%%%%%%%%
%% 4-Results
%%%%%%%%%%%%%%%%%%%%%%%%%%%%%%%%%%%%%%%%%%%%%%%%%%%%%%%%%%%%%%%%%%%%%%%%%%%%%%%%

%%%%%%%%%%%%%%%%%%%%%%%%%%%%%%%%%%%%%%%%%%%%%%%%%%%%%%%%%%%%%%%%%%%%%%%%%%%%%%%%
\section{Numerical Studies}
\label{sec:results}
%%%%%%%%%%%%%%%%%%%%%%%%%%%%%%%%%%%%%%%%%%%%%%%%%%%%%%%%%%%%%%%%%%%%%%%%%%%%%%%%

In this section we analyze our numerical results for the presented model of
\cref{sec:theory} with the adaptive finite element solver
from~\cref{sec:numerics}. Firstly, we introduce and specify the simulation
setup in \cref{subsec:simulation_setup}. Secondly, we consider the numerical
results in detail and discuss the physical effects as well as the numerical
efficiency in \cref{subsec:numerical_results}. For this, we split the analysis
in a 1D spherical symmetric case and a 2D quarter disk of a
nanotube.

%%%%%%%%%%%%%%%%%%%%%%%%%%%%%%%%%%%%%%%%%%%%%%%%%%%%%%%%%%%%%%%%%%%%%%%%%%%%%%%%
\subsection{Simulation Setup}
\label{subsec:simulation_setup}
%%%%%%%%%%%%%%%%%%%%%%%%%%%%%%%%%%%%%%%%%%%%%%%%%%%%%%%%%%%%%%%%%%%%%%%%%%%%%%%%

The derived model in \cref{sec:theory} can be applied to cycle silicon as host
material. The used model parameters as well as the normalized values are listed
in \cref{tab:parameters}. We apply an external lithium flux of $N_\text{ext} =
\SI{1}{\crate}$ for lithiation and $N_\text{ext} = \SI{-1}{\crate}$ for
delithiation.
During the simulations, the particles are cycled between $U_{\max} =
\SI{0.5}{\volt}$ and $U_{\min} = \SI{0.05}{\volt}$ unless otherwise
specified~\cite{kolzenberg2022chemo-mechanical}. This corresponds to an initial
concentration of around $c_0 = 0.02$ and a final time of close to $0.9$, thus
we set $\tfinal = 1.8$ for a total lithiation and delithiation cycle unless
otherwise specified. At the beginning of the delithiation process, what is half
of~$\tfinal$, we continue with our adaptive algorithm, however, we change the
sign of the external lithium flux and enforce for two time steps a time step
size $\tau_n =\num{1e-6}$ as well as an order of one for the temporal
adaptivity. The OCV curve for silicon is chosen
as~\cite{chan2007high-performance}:
\begin{align}
  \label{eq:ocv}
  U_\text{OCV}(z) &=
  \frac{\minus0.2453z^3-0.00527z^2+0.2477z+0.006457}{z+0.002493}.
\end{align}
In the next parts we specify our geometrical reference domain for the
representative battery particle including some further boundary conditions
and symmetry assumption as well as some further implementation details.

\begin{table}[!t]
  \centering
  \caption{Model parameters for numerical experiments
    \cite{kolzenberg2022chemo-mechanical, schoof2022parallelization}.}
  \label{tab:parameters}
  \scalebox{0.9}{
    \begin{tabular}[htbp]{@{}llccc@{}}
      \toprule
      \textbf{Description} & \textbf{Symbol} &
      \multicolumn{1}{c}{\textbf{Value}} & \textbf{Unit} &
      \multicolumn{1}{c}{\textbf{Dimensionless}} \\

      \midrule

      Universal gas constant & $R_\text{gas}$ & $8.314$ &
      \si{\joule\per\mol\per\kelvin}
      & $1$ \\[\defaultaddspace]

      Faraday constant & $F$ & $96485$ & \si{\joule\per\volt\per\mol} & $1$
      \\[\defaultaddspace]

      Operation temperature & $T$ & $298.15$ & \si{\kelvin} & $1$
      \\[\defaultaddspace]

      \midrule

      \multicolumn{5}{c}{Silicon} \\

      \midrule

      Particle length scale & $L_0$ & $\num{50e-9}$ &
      \si{\meter} & 1
      \\[\defaultaddspace]

      Diffusion coefficient & $D$ & $\num{1e-17}$ &
      \si{\square\meter\per\second} &
      $14.4$ \\[\defaultaddspace]

      OCV curve & $U_\text{OCV}$ & $\text{Equation~}\eqref{eq:ocv}$ &
      \si{\volt}
      & $F/R_\text{gas} T \cdot \eqref{eq:ocv}$ \\[\defaultaddspace]

      Young's modulus & $\EH$ & $\num{90.13e9}$ & \si{\pascal} & $116.74$
      \\[\defaultaddspace]

      Poisson's ratio & $\nu$ & $0.22$ & \si{-} & $0.22$ \\[\defaultaddspace]

      Partial molar volume & $v_\text{pmv}$ & $\num{10.96e-6}$ &
      \si{\cubic\meter\per\mol} & $3.41$
      \\[\defaultaddspace]

      Maximal concentration & $\cmax$ & $\num{311.47e3}$ &
      \si{\mol\per\cubic\meter} &
      $1$ \\[\defaultaddspace]

      Initial concentration & $c_0$ & $\num{6.23e3}$ &
      \si{\mol\per\cubic\meter} &
      $\num{2e-2}$
      \\
      \bottomrule
    \end{tabular}
  }
\end{table}

%%%%%%%%%%%%%%%%%%%%%%%%%%%%%%%%%%%%%%%%%%%%%%%%%%%%%%%%%%%%%%%%%%%%%%%%%%%%%%%%
\subsubsection{Geometrical Setup}
\label{subsubsec:geometries}
%%%%%%%%%%%%%%%%%%%%%%%%%%%%%%%%%%%%%%%%%%%%%%%%%%%%%%%%%%%%%%%%%%%%%%%%%%%%%%%%

\begin{figure}[!b]
  \centering
  \includegraphics[width = 0.58\textwidth,
  page=1]{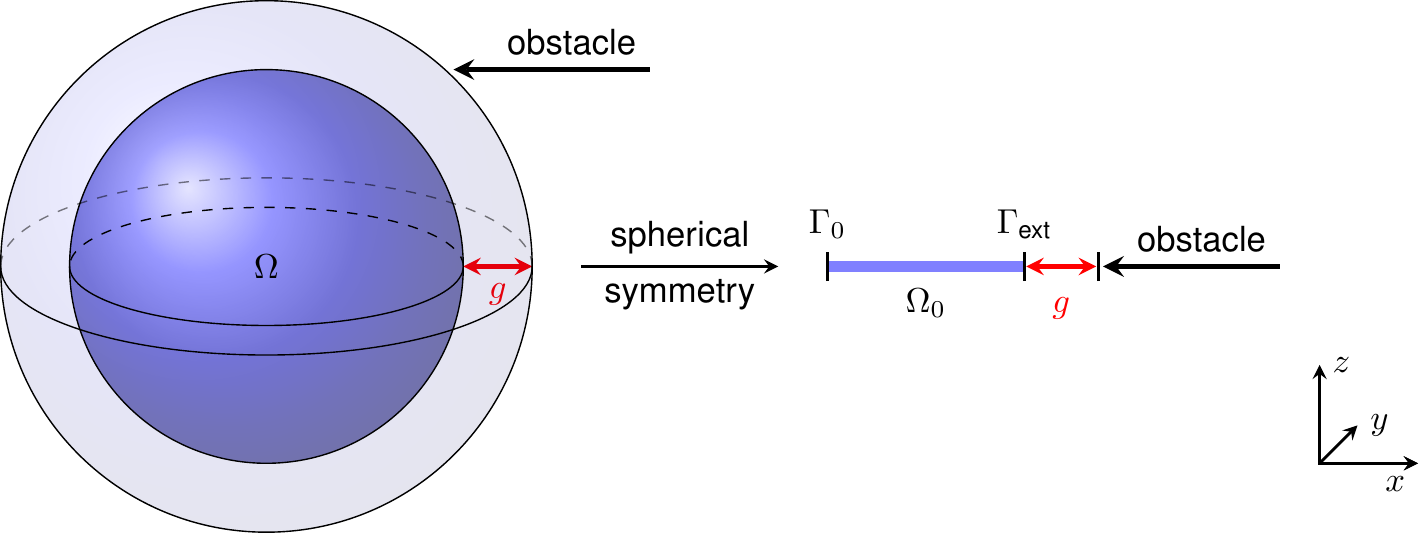}
  \caption{Dimension reduction of a three-dimensional unit sphere with
    surrounded obstacle to the one-dimensional unit interval with spherical
    symmetry and the
    gap function~$g$, based
    on~\cite[Figure~B.1]{castelli2021numerical}.}
  \label{fig:particle_obstacle_1d}
\end{figure}

\begin{figure}[!t]
  \centering
  \includegraphics[width = 0.78\textwidth,
  page=2]{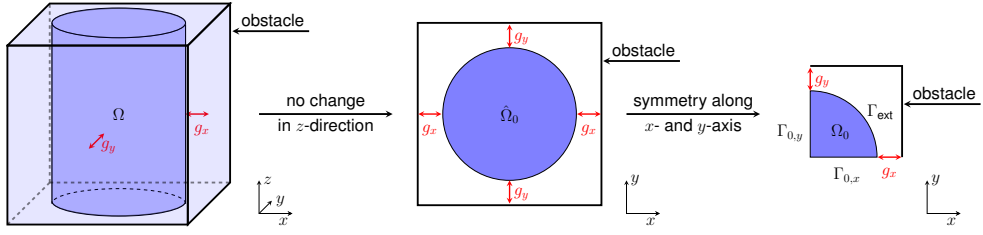}
  \caption{Dimension reduction of a three-dimensional nanowire with surrounded
    rectangular obstacle to the two-dimensional quarter disk
    and the (time-independent) quadratic obstacle
    with the gap function~$(g_x, g_y)^\trp
    = \vect{g} = \hat{\vect{g}} - \vect{X}_0 =  (\hat{g}_x,
    \hat{g}_y)^\trp - (X_{0,x}, X_{0,y})^\trp$ with the components $g_x$
    and
    $g_y$ in $x$- and $y$-directions, respectively.
  }
  \label{fig:particle_obstacle_2d}
\end{figure}

For a representative 3D spherical particle, the computational domain can be
reduced to the 1D unit interval~$\Omega_0 = \left(0,1\right)$ with symmetry
assumptions. The particle is surrounded by an obstacle like a core shell
scenario, compare~\cref{fig:particle_obstacle_1d}. Then, the gap
function~$\vect{g}$
reduces to a one-dimensional parameter $g > 1$. The dimensional reduction
introduces an artificial boundary $\Gamma_0$ at which we impose a no flux
condition for the lithium flux and fixed radial displacement:
\begin{align}
  \vect{N}\cdot \vect{n}_0 = 0, \qquad {u} = {0} \qquad \text{on } (0, \tfinal)
  \times \Gamma_{0}.
\end{align}
At the particle surface $\Gamma_\text{ext}$ the boundary conditions are
considered as
stated in~\cref{sec:theory}. Assuming spherical symmetry we adapt the
quadrature weight to $\de \vect{X}_0 = 4\pi r^2 \de r$ in the discrete finite
element formulation. In this setting, we apply a constant initial concentration
$c_0$, the chemical potential $\mu_0 = \ptl_c \psi_\ch\left(c_0\right)$ and the
one-dimensional stress-free radial displacement ${u}_0(c_0) = r
\left(\lambda_\ch\left(c_0\right) -1 \right)$.

For the 2D simulation we rely on a silicon
nanotube~\cite{chan2007high-performance, wu2019phase}
and reduce the domain to a quarter disk of the nanotube,
see~\cref{fig:particle_obstacle_2d}. Here, we assume symmetry with respect to
the $x$- and $y$-axis and no variations in $z$-direction. The nanotube is
surrounded by a cuboid that will serve as a rigid obstacle. This results in two
artificial boundaries $\Gamma_{0,x}$ and $\Gamma_{0,y}$ in axial direction and
a curved boundary $\Gamma_\text{ext}$ for the surface of the nanotube. On
$\Gamma_\text{ext}$ we apply an isoparametric mapping for the representation
of the curved boundary. No flux conditions and only radial displacement is
expected on the boundaries in axial direction:
\begin{subequations}
  \begin{alignat}{4}
    &\vect{N}\cdot \vect{n}_0 &&= 0, \qquad {u}_y &&= {0},\qquad &&\text{on }
    (0, \tfinal) \times \Gamma_{0,x},\\
    &\vect{N}\cdot \vect{n}_0 &&= 0, \qquad {u}_x &&= {0}, &&\text{on } (0,
    \tfinal) \times \Gamma_{0,y},
  \end{alignat}
\end{subequations}
where $u_i$ is the $i$-th entry of the vector $\vect{u}$, $i\in \{x,y\}$. For
the starting values for the Newton method we choose a constant initial
concentration $c_0$, $\mu_0 = 0$ and $\vect{u}_0 = \vect{0}$.

%%%%%%%%%%%%%%%%%%%%%%%%%%%%%%%%%%%%%%%%%%%%%%%%%%%%%%%%%%%%%%%%%%%%%%%%%%%%%%%%
\subsubsection{Implementation Details}
\label{subsubsec:implementation_details}
%%%%%%%%%%%%%%%%%%%%%%%%%%%%%%%%%%%%%%%%%%%%%%%%%%%%%%%%%%%%%%%%%%%%%%%%%%%%%%%%

For our implementation, we apply an isoparametric fourth-order Lagrangian
finite element method. The basis for our numerical simulation is the finite
element library~\texttt{deal.II}~\cite{arndt2021deal-ii}, implemented
in~\texttt{C++}, together with the interface to the Trilinos
library~\cite[Version~12.8.1]{team2020trilinos} and the UMFPACK
package~\cite[Version~5.7.8]{davis2004algorithm} for the LU-decomposition. All
simulations are executed on a desktop computer with \SI{64}{\giga\byte}~RAM,
Intel~i5-9500~CPU, GCC compiler version 10.3 and the operating system Ubuntu
20.04.5 LTS. The OpenMP Version 4.5 is used for shared memory parallelization
for assembling the Newton matrix, residuals and spatial estimates.

Unless otherwise stated we choose for the space and time adaptive algorithm
the tolerances $\RelTol_t = \RelTol_x = \num{1e-5}$ and $\AbsTol_t = \AbsTol_x
= \num{1e-8}$. For the marking parameters for local coarsening and refinement,
$\theta_{\mathrm{c}} = 0.05$ and $\theta_{\mathrm{r}}=0.5$ are used and a
maximal time step size $\tau_{\max} = \num{1e-2}$.

To get a diagonal structure for the mass matrix $\tens{B}_h$ or $\tens{D}_h$,
respectively, we use mass lumping with a Gau\ss{}--Lobatto quadrature rule,
compare Remark~1 in~\cite{frohne2016efficient}. We initialize the active set
$\activesetk = \emptyset$ and the inactive set $\inactiveset_k = \possibleDOF$
for the lithiation process. During lithiation we check in each time step the
condition~$\hat{{\displacementt}}_p-\hat{{g}}_{p} > 0$ to set this nodal point
active, since $\hat{\lambda}_{p} = 0$ is zero due to the boundary condition for
$p \in \possibleDOF$. During delithiation we check in each time step
$\hat{{\lambda}}_{p} \leq 0$ to set this nodal point inactive again, since
$\hat{{\displacementt}}_p-\hat{{g}}_{p} = 0$. To increase numerical stability
we only change to active points during lithiation and to inactive points during
the following delithiation.

%%%%%%%%%%%%%%%%%%%%%%%%%%%%%%%%%%%%%%%%%%%%%%%%%%%%%%%%%%%%%%%%%%%%%%%%%%%%%%%%
\subsection{Numerical Results}
\label{subsec:numerical_results}
%%%%%%%%%%%%%%%%%%%%%%%%%%%%%%%%%%%%%%%%%%%%%%%%%%%%%%%%%%%%%%%%%%%%%%%%%%%%%%%%
This subsection discusses the numerical simulation results for our two
presented computational domains: a 1D unit interval with modified quadrature
weight to consider a 3D spherical particle and a 2D quarter nanotube. We
analyze the behavior of concentration and stress development inside the
representative active particle and show the efficiency of the adaptive space
and time algorithm for cycling battery active particles with mechanical
constraints.

%%%%%%%%%%%%%%%%%%%%%%%%%%%%%%%%%%%%%%%%%%%%%%%%%%%%%%%%%%%%%%%%%%%%%%%%%%%%%%%%
\subsubsection{1D Spherical Symmetry}
\label{subsubsec:1d_spherical_symmetry}
%%%%%%%%%%%%%%%%%%%%%%%%%%%%%%%%%%%%%%%%%%%%%%%%%%%%%%%%%%%%%%%%%%%%%%%%%%%%%%%%

In this part we consider for one cycle the influence of the obstacle for the
1D spherical symmetric case as shown in \cref{fig:particle_obstacle_1d}.
Firstly, we compare the stress development of a configuration with and without
obstacle. Secondly, have a closer look on the stress development regarding the
radius of the particle as well as the concentration process. Thirdly, we
investigate the influence of the size of the gap function on the stress
development. We
close the part with the consideration of the time step size comparing again a
setup with and without obstacle.

\begin{figure}[t]
  \centering
  \includegraphics[scale=0.65,
  page=1]{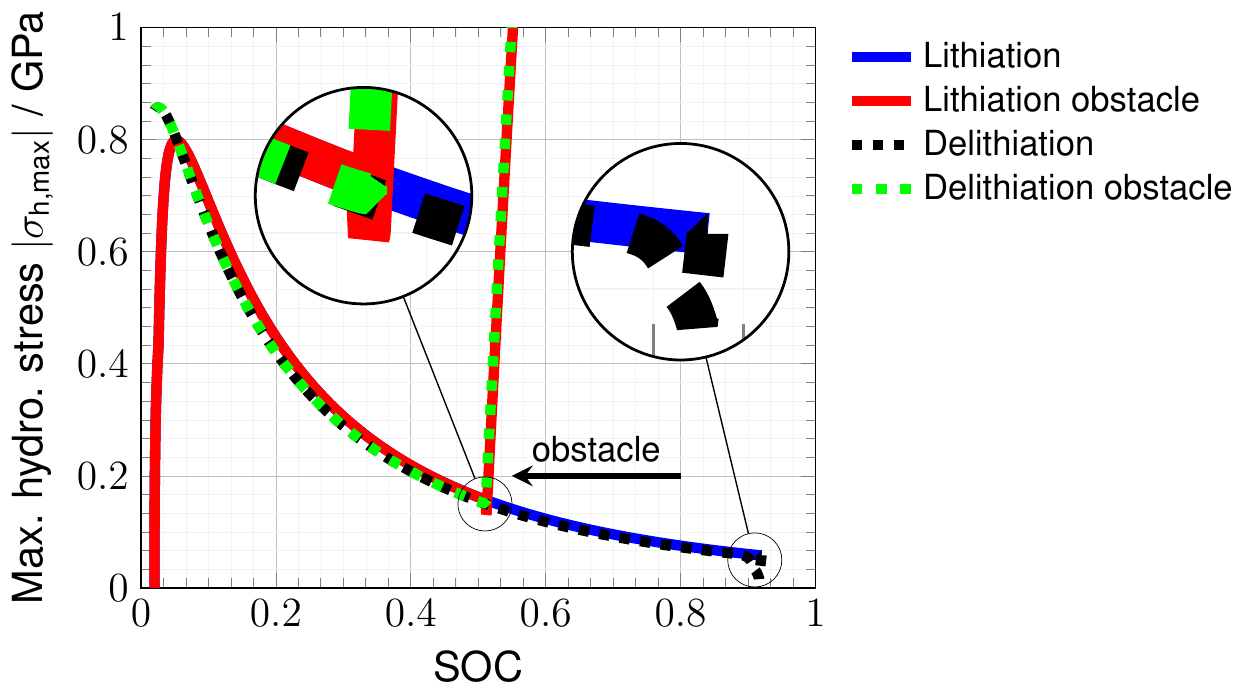}
  \caption{Comparison of the absolute value of the maximal hydrostatic
    stress~$|\sigma_\text{h}|$ over the $\text{SOC}$ for a
    cycle of a 1D spherical symmetric setup without and with obstacle.}
  \label{fig:1d_cycle_max_hydro_sigma}
\end{figure}

\begin{figure}[!b]
  \centering
  \includegraphics[scale=0.6,
  page=1]{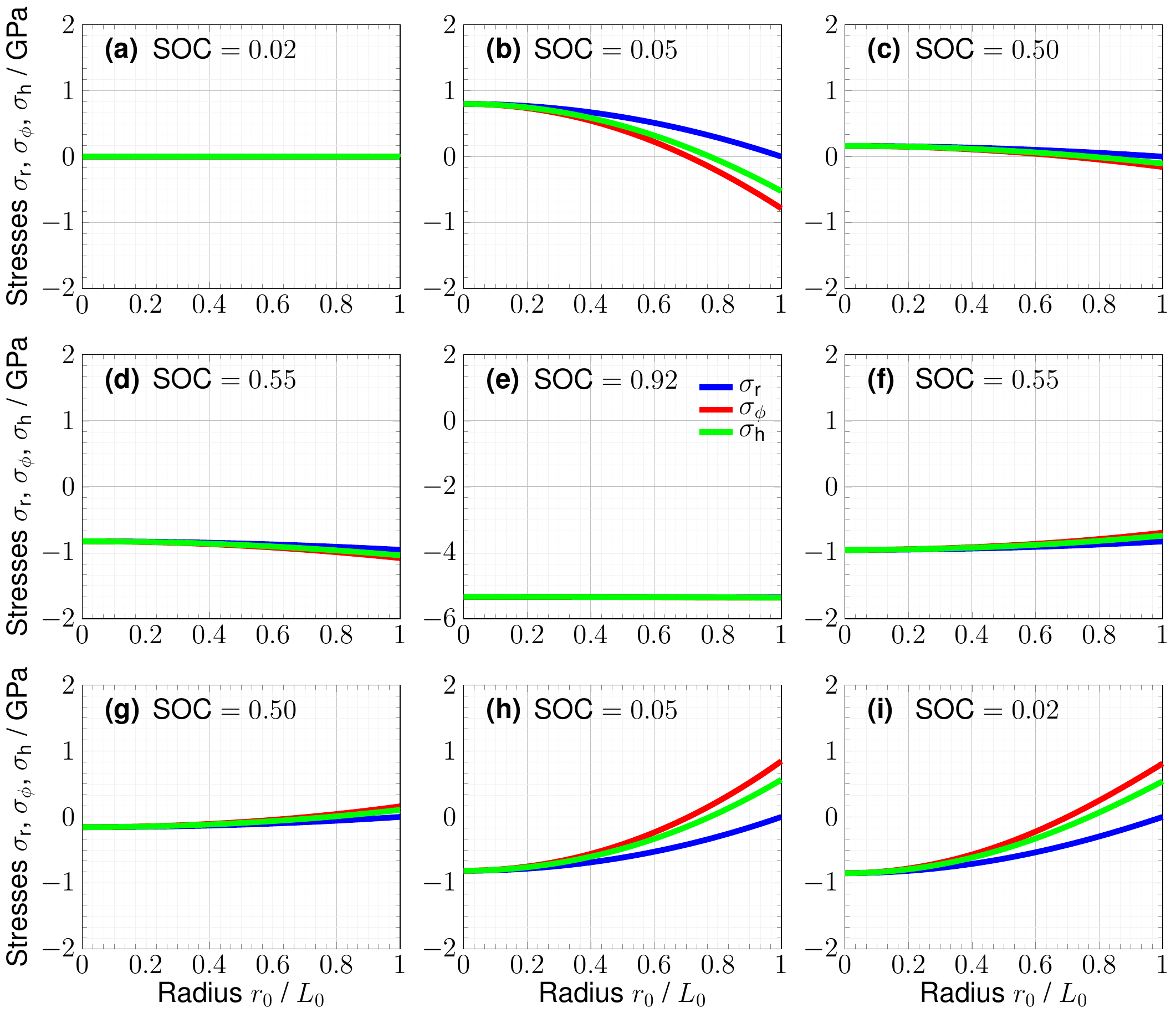}
  \caption{Stress development for the radial $\sigma_\text{r}$ (blue),
    tangential $\sigma_{\phi}$ (red) and hydrostatic part $\sigma_\text{h}$
    (green) over the particle radius for nine different $\text{SOC} \in \{
    0.02,
    0.05, 0.50, 0.55, 0.92, 0.55, 0.50, 0.05, 0.02 \}$ for one cycle
    in~(a)--(i), respectively.}
  \label{fig:1d_cycle_sigma}
\end{figure}

In \cref{fig:1d_cycle_max_hydro_sigma} we see the absolute value of the maximal
hydrostatic stress~$|\sigma_\text{h}| = |1/3 \sigma_\text{r} + 2/3
\sigma_{\phi}|$ in $\si{GPa}$ over the $\text{SOC}$, computed with the radial
and
tangential Cauchy stress, denoted with $\sigma_\text{r}$ and $\sigma_{\phi}$,
respectively. The solid lines represent the lithiation process and the dashed
lines the delithiation process in each case: once without and once with
obstacle. For the gap function, we use~$g=0.4$. Before the
particle gets in contact with the obstacle, the stress curves for the case
without and with obstacle are identical. Shortly after the start of the
lithiation process, a peak of around \SI{0.8}{\giga\pascal} rises. This can be
explained due to the characteristic behavior of the OCV curve, compare Figure~2
in~\cite{schoof2022parallelization}. At around $\text{SOC}=0.51$, the particle
gets
in contact with the obstacle and the stress profiles deviate between the cases.
After a short reduction, the stress values increase significantly until it
reaches a maximum of $|\sigma_{\text{h},\max}| \approx
\SI{5.36}{\giga\pascal}$. The short stress reduction results from the fact,
that all (also positive) stress values inside the particle will change to
compressive stresses which have negative values. So tensile stresses with a
positive value have to go through the zero point. Compare for this also the
detailed analysis in \cref{fig:1d_cycle_sigma}(b) to (c). The stress
development in the case without obstacle flattens out until the change of the
external lithium flux for $\text{SOC}=0.92$ enters. Here again we have a short
drop of the stress value due to the rearrangement from tensile to compressive
values and vice versa. Because of the OCV curve, the stress values increase
again but have a slight shift which can be explained by the rearrangement of
the stress after the change from positive to negative lithium flux. After the
particle detached from the obstacle again, the cases without and with obstacle
coincide again ending at higher stress values.

\begin{figure}[b]
  \centering
  \includegraphics[scale=0.65,
  page=1]{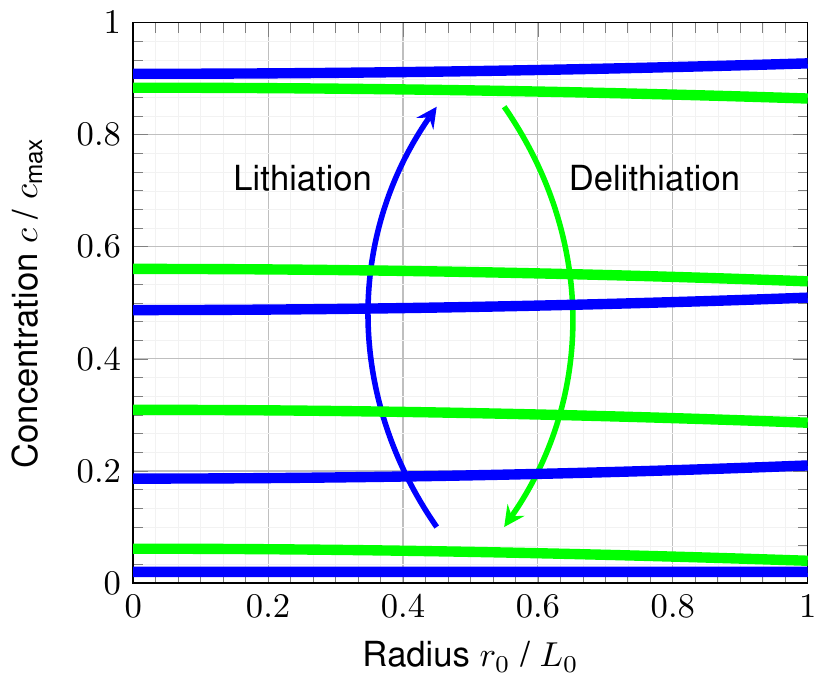}
  \caption{Concentration profile $c$ over the particle radius $r_0$ during one
    charging (blue) and discharging (green) cycle at different
    \mbox{$\text{SOC} \in
      \{0.02, 0.20, 0.50, 0.92, 0.87, 0.55, 0.30, 0.05 \}$}.}
  \label{fig:1d_cycle_concentration}
\end{figure}

In the next step we want to have a closer look on the stress distribution over
the particle radius in the Lagrangian domain. In \cref{fig:1d_cycle_sigma} the
three different stresses, radial, tangential and hydrostatic stress, are
displayed over the particle radius for nine different $\text{SOC}{} \in \{
0.02, 0.05, 0.50, 0.55, 0.92, 0.55, 0.50, 0.05, 0.02 \}$ considering one
cycling. At the initial $\text{SOC}=0.02$ in \cref{fig:1d_cycle_sigma}(a) there
are no stresses present. At $\text{SOC} = 0.05$ the maximal values arise in the
particle center at $r_0=0$. Here, tensile stresses occur whereas at the
particle surface at $r_0=1$ the tangential stresses are compressive. Note the
zero value for the radial stress $\sigma_{\text{r}}$. This is exactly the
stress-free boundary condition in \cref{eq:stress-free_boundary_condition} or
\cref{eq:problem_boundary_constraintse} fulfilled with equality.
In addition, tangential stresses are not equal to zero. For larger
$\text{SOC}$ all stresses decrease due to the influence of the OCV curve.
A good qualitative accordance of our numerical results is given with the
particle stresses in \cite{kolzenberg2022chemo-mechanical} neglecting the SEI
results. At about $\text{SOC}=0.51$ the particle touches the obstacle and the
stress-free condition changes to a Dirichlet boundary condition for the
displacement. The displacement is fixed now, we have
\cref{eq:problem_boundary_constraintsf} with equality. However, now negative
stresses may occur to fulfill \cref{eq:problem_boundary_constraintse} with
strict inequality. Exactly this can be seen in \cref{fig:1d_cycle_sigma}(d)
until the end of the lithiation at $\text{SOC}=0.92$ in
\cref{fig:1d_cycle_sigma}(e).
Here large comprehensive stresses appear throughout the particle domain. Note
the different range on the stress axis in \cref{fig:1d_cycle_sigma}(e).
Changing now the sign of the external flux the discharging begins. In
\cref{fig:1d_cycle_sigma}(f) we are close to the point where the particle
detached from the obstacle again. Compared to \cref{fig:1d_cycle_sigma}(d) the
curvature of the stress profiles is opposite. At $\text{SOC}=0.50$ the particle
has no contact to the obstacle any more and the stress-free boundary condition
comes into effect again. At the end of the delithiation process the stress
values are qualitatively similar to those of the lithiation process but have
the opposite curvature resulting from the negative sign of the external lithium
flux $N_\text{ext}$. In the end in \cref{fig:1d_cycle_sigma}(i) the discharging
process stops at a significant level of stress values compared to the constant
initial concentration. Compared to the stress measurements in
\cite{al-obeidi2016mechanical} investigating coated silicon electrodes, we have
larger stress values. This could be due to the fact that our model is only
based on a chemo-elastic approach. Nevertheless, our model is capable to deal
with the change of the boundary condition during cycling.

\begin{figure}[t]
  \centering
  \includegraphics[scale=0.65,
  page=1]{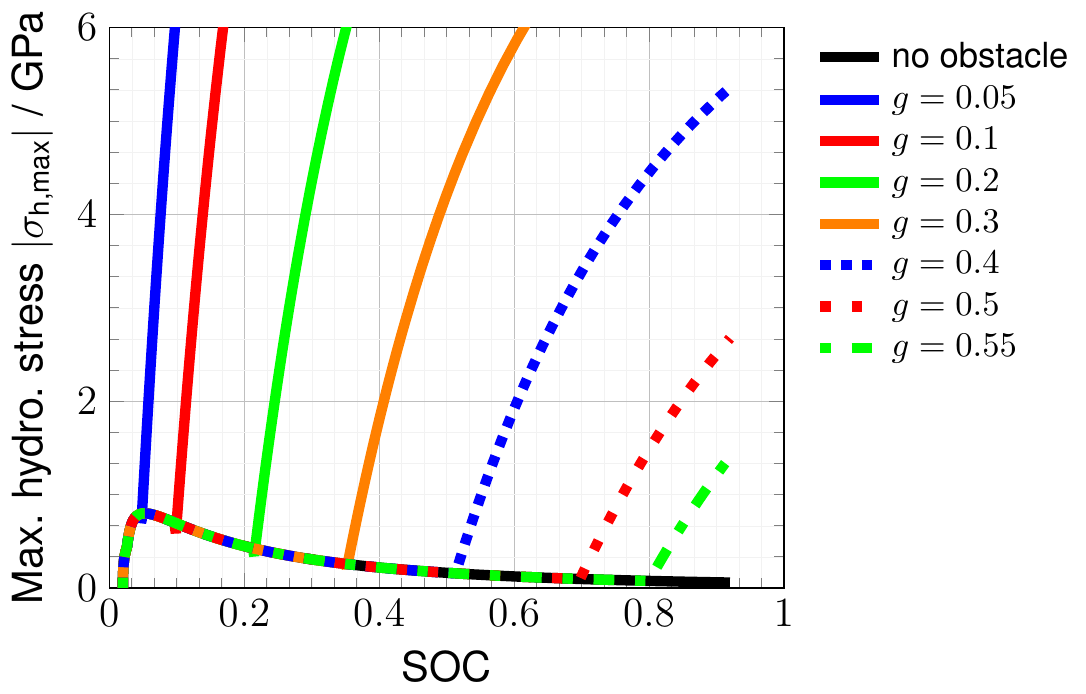}
  \caption{Study on different values for the gap function~$g$: stress
    profiles
    for the
    absolute values of the maximal hydrostatic stress
    $|\sigma_{\text{h},\max}|$
    (colored) over the $\text{SOC}$ compared to the case without obstacle in
    black
    for the charging process only.}
  \label{fig:1d_loop_over_g}
\end{figure}

In \cref{fig:1d_cycle_concentration} the concentration profile at the
\mbox{$\text{SOC} \in \{ 0.02, 0.20, 0.50, 0.92 \}$} for the lithiation in blue
and \linebreak \mbox{$\text{SOC} \in \{ 0.87, 0.55, 0.30, 0.05 \}$} for the
delithiation in green is shown over the particle radius $r_0$. We choose
slightly different $\text{SOC}$ values for the delithiation to have more
distance between the different results. After the constant initial
concentration $c_0$ the concentration profile increases with a slight curvature
until its maximal value at $\text{SOC} = 0.92$. Note that the obstacle contact
has no critical influence on the concentration profile. After switching to the
delithiation process the curvature of the concentration is also opposite like
for the stresses in~\cref{fig:1d_cycle_sigma}. The curvature remains also at
the end of the simulation time for low concentration level.

\cref{fig:1d_loop_over_g} shows the influence on the gap function
on the stress development up to a maximal value of
\SI{6.0}{\giga\pascal} over the $\text{SOC}$. The above investigated case
with~$g=0.4$ is here displayed with the blue dashed line. The simulation
without an obstacle is shown with the solid black line. The smaller the
gap function~$g$
is, the earlier the absolute value of the maximal hydrostatic stress rises.
Likewise, the gradient of the stress profile increases with lower values of the
gap function~$g$. Interestingly, the rise of the stress development is not
constant. For
smaller
gap function values
is the slope significantly higher than
for larger gap function values.
Additionally, the gradient values decrease for higher SOC which could be
explained by the lower curvature of the stress profile itself for larger stress
values, see~\cref{fig:1d_cycle_sigma}(e) compared to,
e.g.,~\cref{fig:1d_cycle_sigma}(c),~(d),~(f)~or~(g).

\begin{figure}[b]
  \centering
  \includegraphics[scale=0.65,
  page=2]{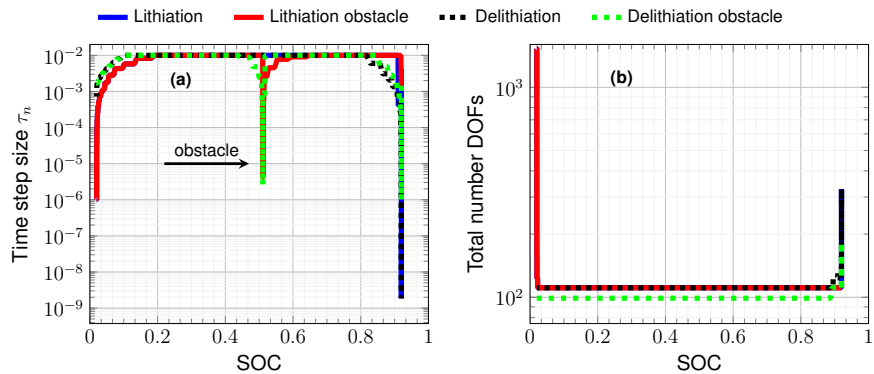}
  \caption{Time step size $\tau_n$ (a) and total number DOFs (b) over
    the $\text{SOC}$ for charging and discharging without and with
    obstacle~($g=0.4$).}
  \label{fig:1d_cycle_time_step_size}
\end{figure}

In the last part of this section, we want to emphasize the need and the
efficiency of the space and time adaptive algorithm. In
\cref{fig:1d_cycle_time_step_size}(a) the time step size $\tau_n$ and in
\cref{fig:1d_cycle_time_step_size}(b) the number of DOFs are plotted over the
$\text{SOC}$ without and with obstacle ($g=0.4$) for the lithiation (solid) and
the delithiation (dashed) process, respectively. We start the simulation with
$\tau_0 =\num{1e-6}$ and $1539$ DOFs in total. After a few time steps the
gradients in concentration, chemical potential and displacement have formed and
the time step size becomes larger until the maximal time step size of
$\tau_{\max} = 0.01$. In the same time, the number of DOFs decreases since no
new gradients occur. However, the time step size decreases significantly over
three orders of magnitudes in the moment when the particle touches the
obstacle, see~\cref{fig:1d_cycle_time_step_size}(a) at $\text{SOC} = 0.51$.
After the particle is in contact with the obstacle, the time step size
increases until $\tau_{\max}$ again. When the delithiation process sets in, we
see three crucial points:
firstly, in \cref{fig:1d_cycle_time_step_size}(a) the time step size of $\tau_n
= \num{1e-6}$ for the delithiation with obstacle seems to be enough. Secondly,
this is in contrast to the case without obstacle. Here, a drop of order of
magnitude over more than six is needed compared to the maximal time step size
$\tau_{max}$. Thirdly, the number of DOFs behave in a similar way. More DOFs
are needed in the case without obstacle compared to the case with
obstacle. An explanation might be the again the lower curvature in the obstacle
case so the changes in the gradients are not so large compared to the
obstacle-free situation. After all gradients feature a reversed direction,
the time step size $\tau_n$ increases as well as the number of DOFs
decreases again. In the obstacle case the number of DOFs is even slightly lower
compared to the simulation without obstacle. However, the time step
size~$\tau_n$ drops down to approximately the same level as for the charging
process when the particle detaches from the obstacle. Finally, the time step
size $\tau_n$ flattens out at the end of the simulation time. All in all, the
efficiency of the spatial and temporal adaptive algorithm is clearly visible
and of significant importance due to the switching point of the obstacle
contact and the switching point of the lithium flux. Without the adaptivity in
space and time we would have to use the lowest time step size and the highest
number of DOFs throughout the total simulation. See for more details about the
numerical efficiency for phase-field
materials~\cite[Section~4.3]{castelli2021efficient}.

%%%%%%%%%%%%%%%%%%%%%%%%%%%%%%%%%%%%%%%%%%%%%%%%%%%%%%%%%%%%%%%%%%%%%%%%%%%%%%%%
\subsubsection{2D Quarter Nanotube}
\label{subsubsec:2d_quarter_nanotube}
%%%%%%%%%%%%%%%%%%%%%%%%%%%%%%%%%%%%%%%%%%%%%%%%%%%%%%%%%%%%%%%%%%%%%%%%%%%%%%%%

Here we analyze the numerical results of the 2D quarter disk
as described in
\cref{subsec:simulation_setup}. For this simulation we use the parameters
$\theta_{\mathrm{c}} = 0.005$, $\RelTol_t = \RelTol_x = \num{4e-5}$, $\AbsTol_t
= \AbsTol_x = \num{4e-8}$. Furthermore, we choose $\tau_0 = \num{1e-8}$,
$\tau_{\max} = \num{1e-3}$ and $\tfinal = 0.2$ to get an appropriate cycling
and use a constant grid and a constant time step size for two time steps after
the discharge process is started. To increase numerical stability during the
time steps with active points, we accept this time step after one spatial
refinement if the Newton update criterion is fulfilled and also in one
following time step, we allow to skip the spatial refinement criterion.
The time-independent~$\hat{\vect{g}}$
is defined
by~$\hat{\vect{g}} = (\hat{g}_x, \hat{g}_y)^\trp = (1.07, 1.07)^\trp$.

\begin{figure}[!b]
  %\centering
  \hspace{-0.375cm}
  \includegraphics[scale=0.59,
  page=1]{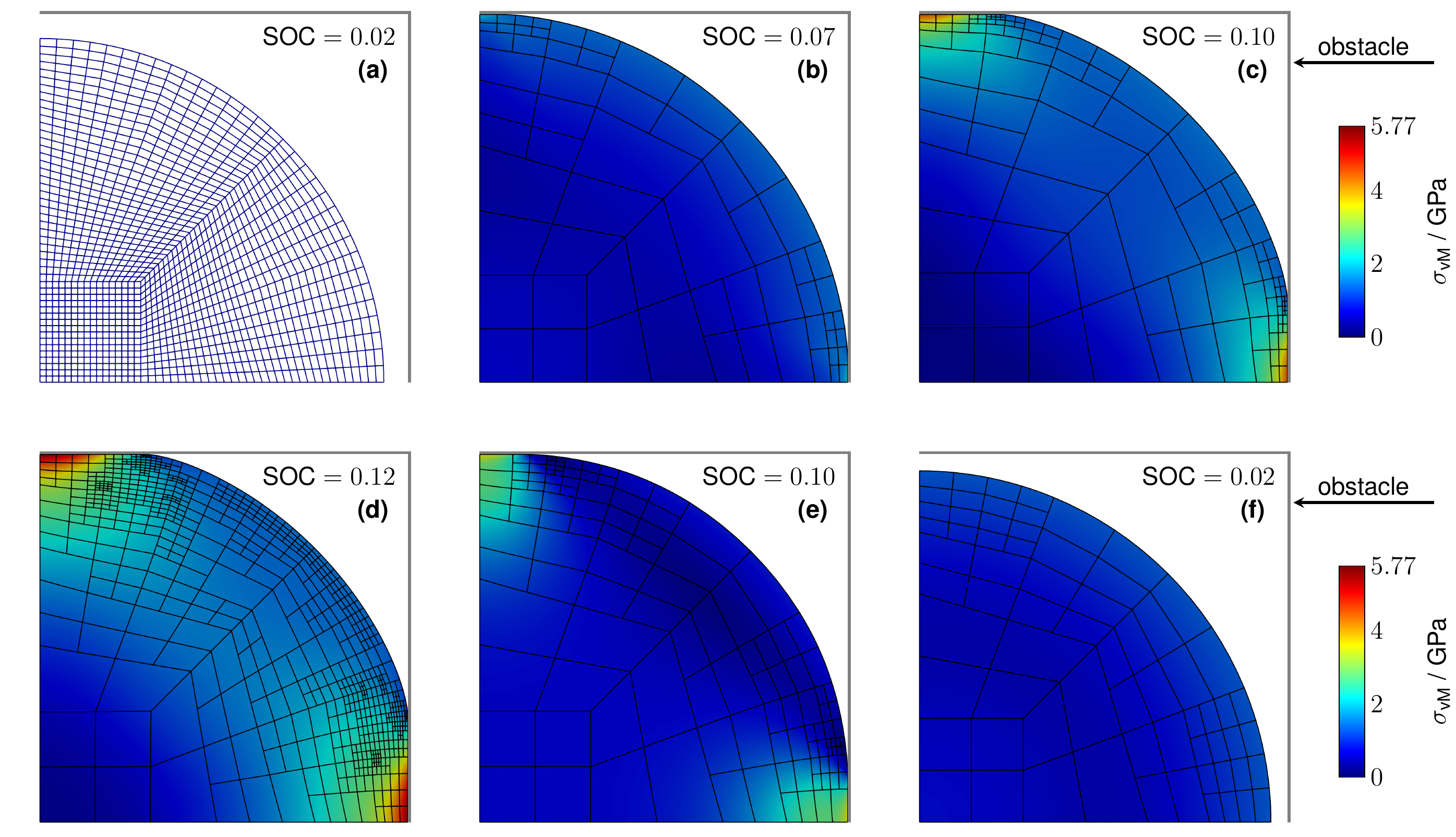}
  \caption{Von Mises stress $\sigma_{\text{vM}}$ of the two-dimensional quarter
    disk of a nanotube in the Eulerian domain $\Omega$
    surrounded by a square
    shaped obstacle at different
    $\text{SOC} \in \{0.02, 0.07, 0.10, 0.12, 0.10, 0.02\}$ for charging and
    discharging with $\tfinal = 0.2$.}
  \label{fig:2d_cycle_sigma}
\end{figure}

\begin{figure}[!t]
  \centering
  \includegraphics[scale=0.65,
  page=1]{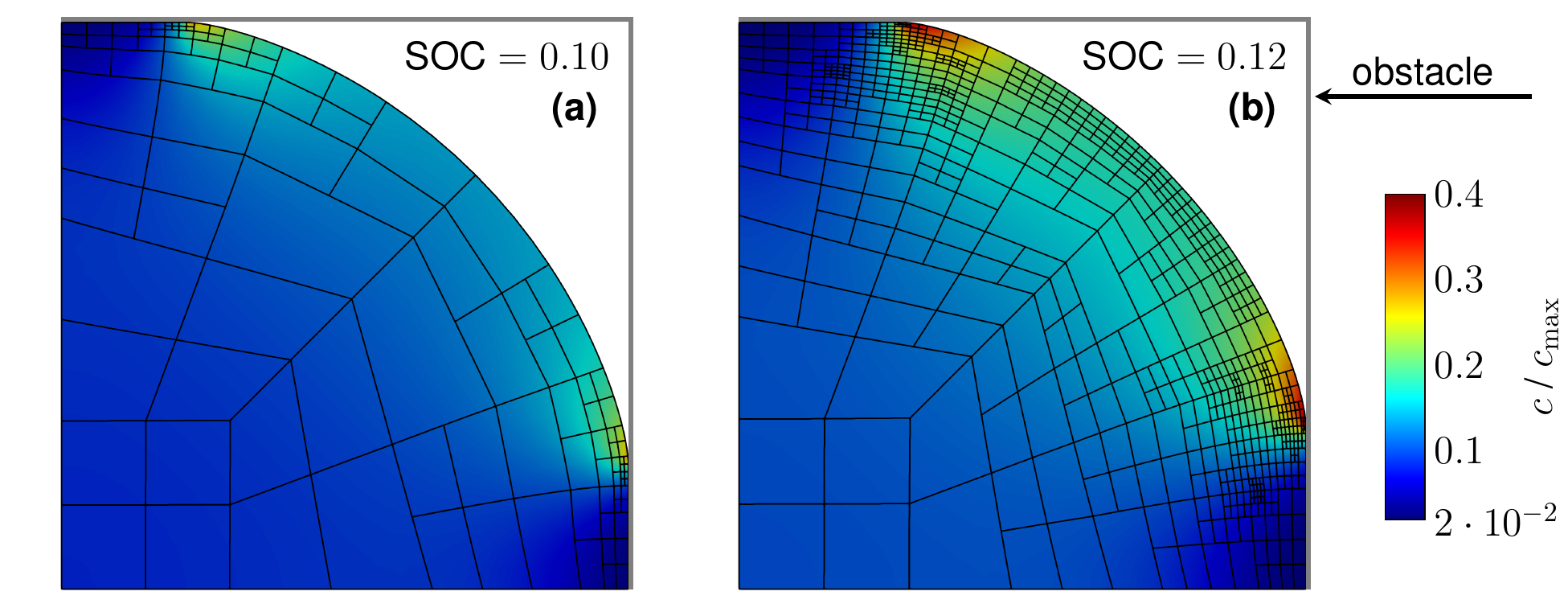}
  \caption{Concentration $c$ of the two-dimensional quarter
    disk of a
    nanotube in the Eulerian domain $\Omega$ surrounded by a square shaped
    obstacle at~$\text{SOC} = 0.10, 0.12$.}
  \label{fig:2d_concentration}
\end{figure}

In \cref{fig:2d_cycle_sigma}(a)--(f) the von Mises stress in the general plane
state
\begin{align}
  \label{eq:von_mises_stress}
  \sigma_{\text{vM}} = \sqrt{\sigma_{11}^2 + \sigma_{22}^2 -
    \sigma_{11}\sigma_{22} + 3\sigma_{12}^2}
\end{align}
is displayed for six different $\text{SOC} \in \{0.02, 0.07, 0.10, 0.12, 0.10,
0.02\}$
warped by the displacement vector to the Eulerian domain $\Omega$ and
surrounded by the obstacle. The uniform grid for the initial
time step without any stresses is shown in~\cref{fig:2d_cycle_sigma}(a). At
$\text{SOC} = 0.07$ there are twelve active points, six at the lower right
corner and six at the upper left corner. At this state all stresses are below
\SI{2.0}{\giga\pascal} but it is visible that the largest stresses occur at the
contact points. This observation strengthens for higher~$\text{SOC}$. The
highest stress values occur at the first contact points, compare, e.g.,
\cref{fig:2d_cycle_sigma}(c) or~(d). This results from the suppression of
the volume increase of the host material. Note that we charge and discharge
with a constant lithium flux $N_\text{ext}$. At this point a Butler--Volmer
boundary condition might be more appropriate but we postpone this to future
work. Near new contact points we have a higher grid resolution due to changes
mainly in the concentration profile. This point is discussed in more detail in
\cref{fig:2d_concentration}. \cref{fig:2d_cycle_sigma}(d) is at~$\text{SOC} =
0.119995$, shortly after the discharging process was started. Here, the
largest stress values occur which are again larger than the measured stresses
in \cite{al-obeidi2016mechanical}. Similar to the 1D simulation a higher grid
resolution and small time steps appear due to the change in the sign of the
constant lithium flux $N_{\text{ext}}$. Again, the efficiency of the space and
time adaptivity is crucial to appropriately capture the change in the physics.
At $\text{SOC}=0.10$ of the delithiation process in \cref{fig:2d_cycle_sigma}(e)
the grid has a coarser structure again. However, the stress distribution differs
from $\text{SOC}=0.10$ of the lithiation process. The occurring maximal stress
values are lower and the distribution of the high values is more orientated
towards the axes-direction instead of the obstacle direction. Moreover, the
particle sections which are in contact with the obstacle are smaller, too. We
also take notice of a dark blue region of low stresses that appears orthogonal
to the first bisector of the coordinate system. In \cref{fig:2d_cycle_sigma}(f)
the final time $\tfinal=0.2$ is reached with state of low stresses.

\begin{figure}[!b]
  \centering
  \includegraphics[scale=0.65,
  page=1]{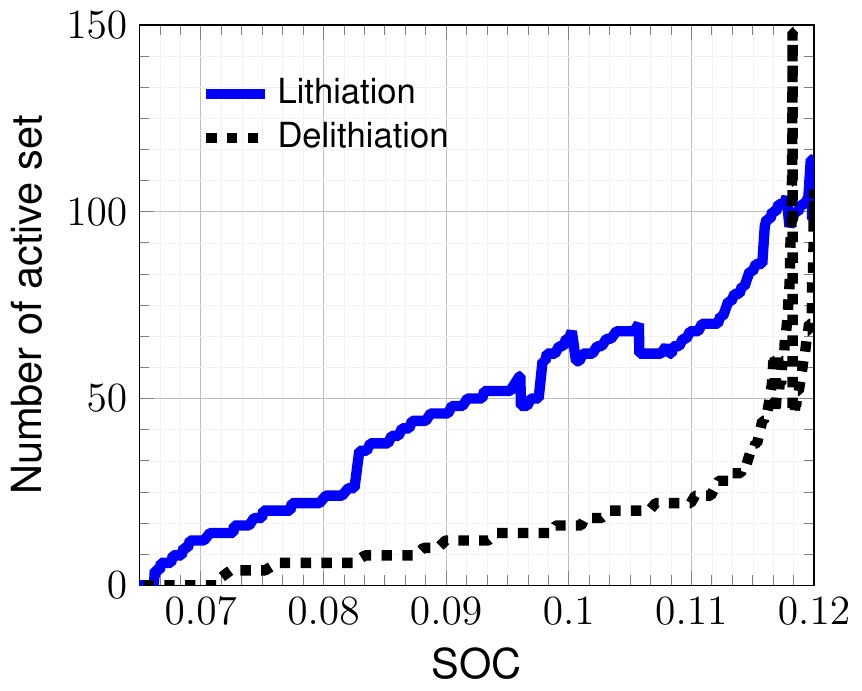}
  \caption{Number of active set $\activeset$ over the part of the $\text{SOC}$
    with a
    nonempty active set of the two-dimensional quarter disk
    of a nanotube
    surrounded by a square shaped obstacle for charging (solid blue) and
    discharging (dashed black) with $\tfinal = 0.2$.}
  \label{fig:2d_number_active_set}
\end{figure}

The concentration $c$ for the two $\text{SOC} $-values $0.10$ and $0.12$ is
shown in \cref{fig:2d_concentration}. The result
in~\cref{fig:2d_concentration}(a) emerged by the charging case, whereas
\cref{fig:2d_concentration}(b) arose at the same $\text{SOC} = 0.119995$ as
in~\cref{fig:2d_cycle_sigma}(d). Eye-catching is the concentration profile near
the obstacle contact. The region where the particle is in contact with the
obstacle has a significantly lower concentration value compared with the region
where the particle in not yet in touch with the obstacle. This effect is even
more pronounced in~\cref{fig:2d_concentration}(b) than
in~\cref{fig:2d_concentration}(a).
Compared to the 1D spherical symmetry setup we have now a different behavior
for the concentration around the obstacle. The region contacting the obstacle
has lower concentration values possibly because the free energy density is
smaller in the contact region with a chemical part with lower concentration
values and a larger elastic part due to the obstacle contact. Also the gradient
in the concentration profile as well as the one in the chemical potential has
to be reverted leading to the finer grid distribution during the discharging
process.

\cref{fig:2d_number_active_set} presents the number of DOFs of the active set
$\activeset$ over the $\text{SOC}$ when the particle is in contact with the
obstacle during cycling. A hysteresis in the number of DOFs is clearly visible.
Specifically, the first active points arise at an earlier $\text{SOC}$ whereas
the last active points vanish at a larger $\text{SOC}$. The approximately linear
increase of the number of active points is interrupted in the middle part by
some sudden peaks and lows due to the refining and coarsening mechanism. In the
end of charging the increase is again nearly linear except for two refinement
peaks. During delithiation the number of active points has one large peak
at~$\text{SOC} =0.118252$ resulting from the change of the gradient of the
concentration and the chemical potential. This change is located a little
further inside the particle and not directly at the boundary and therefore has
a little delay in time. After the curvatures are rearranged the number of
active points decrease with a smaller slope compared to the charging
until the particle detaches from the obstacle. The total number of DOFs varies
between a range of 3596 and 96164 DOFs and time step sizes $\tau_n$ between
$\num{1.58e-8}$ and $\tau_{\max} = \num{1e-3}$. The total computational
time of this two dimensional setup is less than $37$ minutes using the strength
of the semismooth space and time adaptive algorithm.

%%%%%%%%%%%%%%%%%%%%%%%%%%%%%%%%%%%%%%%%%%%%%%%%%%%%%%%%%%%%%%%%%%%%%%%%%%%%%%%%
%% End of 4-results.tex
%%%%%%%%%%%%%%%%%%%%%%%%%%%%%%%%%%%%%%%%%%%%%%%%%%%%%%%%%%%%%%%%%%%%%%%%%%%%%%%%

  %% Summary and conclusion
  % !TeX encoding = UTF-8
% !TeX spellcheck = en_US

%%%%%%%%%%%%%%%%%%%%%%%%%%%%%%%%%%%%%%%%%%%%%%%%%%%%%%%%%%%%%%%%%%%%%%%%%%%%%%%%
%% 5-Conclusion
%%%%%%%%%%%%%%%%%%%%%%%%%%%%%%%%%%%%%%%%%%%%%%%%%%%%%%%%%%%%%%%%%%%%%%%%%%%%%%%%

\section{Summary and Conclusion}
\label{sec:conclusion}

%% Summary
We have developed a thermodynamically consistent chemo-mechanical model for
battery active particles coupling chemical effects in the energy density
dependent on a measured OCV curve and finite deformations together with
mechanical boundary constraints of the obstacle problem during cycling of
lithium-ion batteries. Furthermore, we have combined the primal-dual active set
strategy as semismooth Newton method to the space and time adaptive solution
algorithm with higher-order finite elements for the numerical simulation of our
model equations. Using silicon as an example for a host material we have
investigated several simulation setups in one and two dimensions to discuss the
simulation results from a physical and numerical point of view.
We have figured out that the stresses increase significantly if the particle
has only limited surrounding space and is in contact with the obstacle. The
distance to the obstacle has a crucial influence on the slope of the stress
increase during cycling. Because of the switch in the sign of the external
lithium flux the curvature of the concentration, chemical potential and stress
profiles also have to rearrange oppositely resulting in a hysteresis
development of the concentration, chemical potential and stress profiles.
Although in the one-dimensional simulation setup the obstacle has almost no
influence on the concentration over the particle radius, a clear difference can
be seen in the two-dimensional case: a new lithium-poor region occurs near the
obstacle area reducing the energy density due to the large ratio of the elastic
part. In the two-dimensional case, the largest stress values occur near the
first contact area between the active particle and the obstacle and a clear
hysteresis of the stress values can be seen.

%% Results numerical
Looking at the time step scale the power of the adaptive method is revealed
immediately. Without a variable time step size and order, the simulation must
run with the smallest present time step size, e.g., using a standard backward
Euler scheme to correctly capture all physical effects. This would result
in a significant increase of computational costs compared to our numerical
solution procedure. This is especially crucial for the change of the sign for
the lithium flux to simulate a total cycling and also for long term battery
operations. Similarly, the spatial adaptivity is necessary to appropriately
capture all physical mechanisms, especially the new phenomena of the
lithium-poor phase around the obstacle contact in the two-dimensional setup.
The usage of the semismooth Newton method does not increase the number of DOFs
of the linear system solving for the Newton update and makes it very useful for
higher-dimensional computations~\cite{brunssen2007fast}.

%% Outlook
The efficient extension to various two- and three-dimensional geometries with
differently shaped particles and surrounded obstacles can be investigated in
future works together with long term battery cycles. The large emerging
stresses can lead to the need of further coupling, e.g., with plastic effects
or additional fracture mechanisms. All this together can help in the
understanding of mechanical degradation, capacity fade and battery aging.
The investigation of numerically expensive simulations like for phase
separation materials as LFP or LMO with surrounded obstacles is also another
promising application for this highly efficient adaptive solver.

%%%%%%%%%%%%%%%%%%%%%%%%%%%%%%%%%%%%%%%%%%%%%%%%%%%%%%%%%%%%%%%%%%%%%%%%%%%%%%%%
%% End of 5-conclusion.tex
%%%%%%%%%%%%%%%%%%%%%%%%%%%%%%%%%%%%%%%%%%%%%%%%%%%%%%%%%%%%%%%%%%%%%%%%%%%%%%%%

  %%%%%%%%%%%%%%%%%%%%%%%%%%%%%%%%%%%%%%%%%%%%%%%%%%%%%%%%%%%%%%%%%%%%%%%%%%%%%%
  %% Declaration of competing interest
  %%%%%%%%%%%%%%%%%%%%%%%%%%%%%%%%%%%%%%%%%%%%%%%%%%%%%%%%%%%%%%%%%%%%%%%%%%%%%%
  % The usual declaration
  \section*{Declaration of competing interest}

  \noindent
  The authors declare that they have no known competing financial interests or
  personal relationships that could have appeared to influence the work in this
  paper.

  %%%%%%%%%%%%%%%%%%%%%%%%%%%%%%%%%%%%%%%%%%%%%%%%%%%%%%%%%%%%%%%%%%%%%%%%%%%%%%
  %% CRediT authorship contribution statement
  %%%%%%%%%%%%%%%%%%%%%%%%%%%%%%%%%%%%%%%%%%%%%%%%%%%%%%%%%%%%%%%%%%%%%%%%%%%%%%
  % CRediT authorship contribution statement, see
  %
  %https://www.elsevier.com/authors/policies-and-guidelines/credit-author-statement
  \section*{CRediT authorship contribution statement}

  \noindent
  \textbf{R. Schoof:}
  Methodology, Software, Validation, Formal analysis, Investigation,
  Data Curation, Writing -- original draft, Visualization.
  \textbf{G. F. Castelli:}
  Software, Data Curation, Writing -- review \& editing.
  \textbf{W. Dörfler:}
  Conceptualization, Resources, Writing -- review \& editing, Supervision,
  Project administration, Funding acquisition.

  %%%%%%%%%%%%%%%%%%%%%%%%%%%%%%%%%%%%%%%%%%%%%%%%%%%%%%%%%%%%%%%%%%%%%%%%%%%%%%
  %% Acknowledgment
  %%%%%%%%%%%%%%%%%%%%%%%%%%%%%%%%%%%%%%%%%%%%%%%%%%%%%%%%%%%%%%%%%%%%%%%%%%%%%%
  \section*{Acknowledgement}

  \noindent
  The authors thank L. von Kolzenberg and L. Köbbing
  for intensive and constructive discussions about modeling silicon particles.
  R.S. and G.F.C. acknowledge financial support by the German Research
  Foundation~(DFG) through the Research Training Group 2218
  SiMET~--~Simulation of Mechano-Electro-Thermal processes in Lithium-ion
  Batteries, project number 281041241.

  %%%%%%%%%%%%%%%%%%%%%%%%%%%%%%%%%%%%%%%%%%%%%%%%%%%%%%%%%%%%%%%%%%%%%%%%%%%%%%
  %% ORCID
  %%%%%%%%%%%%%%%%%%%%%%%%%%%%%%%%%%%%%%%%%%%%%%%%%%%%%%%%%%%%%%%%%%%%%%%%%%%%%%
  \section*{ORCID}

  \noindent
  R. Schoof: \url{https://orcid.org/0000-0001-6848-3844}

  \noindent
  G. F. Castelli: \url{https://orcid.org/0000-0001-5484-6093}

  \noindent
  W. D\"orfler: \url{https://orcid.org/0000-0003-1558-9236}

  %%%%%%%%%%%%%%%%%%%%%%%%%%%%%%%%%%%%%%%%%%%%%%%%%%%%%%%%%%%%%%%%%%%%%%%%%%%%%%
  %% Appendix
  %%%%%%%%%%%%%%%%%%%%%%%%%%%%%%%%%%%%%%%%%%%%%%%%%%%%%%%%%%%%%%%%%%%%%%%%%%%%%%

  %% Appendix: Abbreviations, Symbols, Tensor Analysis
  % !TeX encoding = UTF-8
% !TeX spellcheck = en_US

%%%%%%%%%%%%%%%%%%%%%%%%%%%%%%%%%%%%%%%%%%%%%%%%%%%%%%%%%%%%%%%%%%%%%%%%%%%%%%%%
%% Appendix
%%%%%%%%%%%%%%%%%%%%%%%%%%%%%%%%%%%%%%%%%%%%%%%%%%%%%%%%%%%%%%%%%%%%%%%%%%%%%%%%

\appendix

\addcontentsline{toc}{section}{Appendices}
\section*{Appendices}

%%%%%%%%%%%%%%%%%%%%%%%%%%%%%%%%%%%%%%%%%%%%%%%%%%%%%%%%%%%%%%%%%%%%%%%%%%%%%%%%
\renewcommand{\thesection}{A}
\section{Abbreviations and Symbols}
\label{app:abbreviations_and_symbols}
%%%%%%%%%%%%%%%%%%%%%%%%%%%%%%%%%%%%%%%%%%%%%%%%%%%%%%%%%%%%%%%%%%%%%%%%%%%%%%%%
\hspace*{-0.60cm}
\begin{longtable}[l]{p{2.0cm}p{5.5cm}p{2.5cm}p{5.5cm}}
  \multicolumn{4}{l}{\textbf{Abbreviations}} \\
  \ \\
  DAE & differential algebraic equation &
  NCP & nonlinear complementary problem \\
  DOF & degree of freedom &
  OCV & open circuit voltage \\
  KKT & Karush--Kuhn--Tucker &
  SOC & state of charge
\end{longtable}

\hspace*{-0.60cm}
\begin{longtable}[l]{p{2.0cm}p{5.5cm}p{2.0cm}p{5.5cm}}
  \textbf{Symbol} & \textbf{Description}	& & \\
  \ \\
  \multicolumn{2}{l}{Latin symbols} & \multicolumn{2}{l}{Greek symbols}
  \vspace{0.2cm}\\
  ${}_{z}\tens{A}_h(z, \tens{G})$ & partial derivative of $\tens{A}_h$
  regarding $z$ &
  $\alpha > 0$ & coefficient in NCP function \\
  ${}_{\tens{G}}\tens{A}_h(z, \tens{G})$ & partial derivative of $\tens{A}_h$
  regarding $\tens{G}$ &
  $\alpha_{\tau_n} > 0$ & coefficient for adaptive time\\
  $\activeset$ & set of all active nodes of $\possibleDOF$&
  & discretization \\
  $\tens{B}_h$ & discrete auxiliary matrix &
  $\delta_{j,k}$ & Kronecker delta \\
  $\mathcalC$ & NCP function
  &
  $\Gamma_\possibleDOF$ & potential contact zone \\
  $\Ctensor$ & fourth-order stiffness tensor &
  $\vect{\lambda}$ & Lagrange multiplier \\
  $\tens{D}_h$ & discrete diagonal matrix &
  $\lambda_\ch$ & factor of concentration induced \\
  $\tens{E}_\el$ & elastic strain tensor &
  & deformation gradient  \\
  $\tens{F} = \grad_0\vect{\chi}$ & deformation gradient tensor &
  $\mu$ & chemical potential \\
  $\tens{F} = \tens{F}_\ch\tens{F}_\el$ & multiplicative decomposition of
  $\tens{F}$ &
  $\nu$ & Poisson's ratio \\
  $\tens{F}_\ch$ & chemical deformation gradient &
  $\Omega$ & Eulerian domain \\
  $\tens{F}_\el$ & elastic deformation gradient &
  $\Omega_0$ & Lagrangian domain \\
  $\vect{g}$ & gap function &
  $\varphi$ & scalar valued test function \\
  $\hat{\vect{g}}$ & projection &
  $\psi$ & total free energy density \\
  $\inactiveset$ & set of all inactive nodes of $\possibleDOF$ &
  $\psi_\ch$ & chemical part of free energy density \\
  ${m}$ & scalar valued mobility &
  $\psi_\el$ & elastic part of free energy density \\
  $\vect{n}$, $\vect{n}_0$ & normal vector on $\Omega$, $\Omega_0$ &
  $\boldsymbol{\sigma}$ & Cauchy stress tensor\\
  $N$, $N_\Lambda$ & number of nodes of ${V}_h$, $\vect{\Lambda}_h$ &
  $\vect{\xi}_j$ & vector valued test function of node $j$\\
  $\vect{N}$ & lithium flux &
  ${\xi}_j$ & scalar basis function: nonzero entry \\
  ${N}_\text{ext}$ & external lithium flux &
  & of $\vect{{\xi}_j}$ of node $j$ 	\\
  $\otherDOF$ & set of all other nodes of $\OmegaL$: $\alldofs \textbackslash
  \possibleDOF$ &
  \multicolumn{2}{l}{Mathematical symbols}	\\
  $\tens{P}$ & first Piola--Kirchhoff stress tensor &
  $\partial \Omega_0$ & boundary of $\Omega_0$  \\
  $\possibleDOF$ & set of all potential contact nodes on
  $\Gamma_\possibleDOF$ &
  $\grad_0$ & gradient vector in Lagrangian domain \\
  $\alldofs$ & set of all nodes on $\Omega_0$ &
  $\square \!:\! \tilde{\square}$ & reduction of two dimensions of \\
  $U_\text{OCV}$ & OCV curve &
  & two tensors $\square$ and $\tilde{\square}$\\
  $\displacement = \vect{x} - \vect{X}_0$ & displacement vector &
  \multicolumn{2}{l}{Indices} \\
  $\displacement_h$ & discrete displacement vector or &
  ${\square}_{0}$ & considering variable in Lagrangian \\
  & algebraic representation &
  & domain or initial time \\
  ${V}$ & scalar valued function space &
  ${\square}_{h}$ & finite dimensional function of $\square$ or \\
  $\vect{V}$ & vector valued function space &
  & algebraic representation of $\square$ with \\
  $\vect{V}^*$ & subset of $\vect{V}$ &
  & respect to basis function \\
  $\vect{V}^+$ &  subset of $\vect{V}^*$ &
  ${\square}_{\possibleDOF}$ & discrete vector with all entries on set
  $\possibleDOF$\\
  $\vect{x} = \vect{\chi}\left(t, \vect{X_0}\right) $ & motion &
  ${\partial}_{\square}$ & partial derivative with respect to $\square$ \\
  $\vect{X}_0$ & space coordinate in Lagrangian domain &
  & \\
\end{longtable}

%%%%%%%%%%%%%%%%%%%%%%%%%%%%%%%%%%%%%%%%%%%%%%%%%%%%%%%%%%%%%%%%%%%%%%%%%%%%%%%%
\renewcommand{\thesection}{B}
\section{Tensor Analysis}
\label{app:tensor_analysis}
%%%%%%%%%%%%%%%%%%%%%%%%%%%%%%%%%%%%%%%%%%%%%%%%%%%%%%%%%%%%%%%%%%%%%%%%%%%%%%%%
We use the following notation for a scalar $a \in \IR$, first-order vectors
$\vect{b}$, $\vect{c}\in \IR^d$, second-order tensors $\tens{A}$, $\tens{B}$,
$\tens{C} \in \IR^{d,d}$, the second-order identity $\tens{Id} \in
\IR^{d,d}$ and a fourth-order tensor~$\Ctensor \in \IR^{d,d,d,d}$:
\begin{align}
  a &= \vect{b} \cdot \vect{c}, \qquad \tens{A} = \tens{B}\tens{C}, \qquad
  \tens{A} = \Ctensor\left[ \tens{B}\right],
  \qquad a = \tens{A} \!:\! \Ctensor \left[ \tens{B}
  \right]\\
  a &
  = \tens{B} \!:\! \tens{C}
  = \text{tr}\left(\tens{B}^\trp \tens{C}\right)
  = \text{tr}\left(\tens{C}^\trp \tens{B}\right)
  =\text{tr}\left(\tens{B} \tens{C}^\trp\right)
  =\text{tr}\left(\tens{C}^\trp \tens{B}\right)
  = \tens{C} \!:\! \tens{B},
\end{align}
where $\text{tr}\left(\tens{A}\right)$ denotes the trace of a tensor~$\tens{A}$.
Further, we write for first-order vectors $\vect{a}$, $\vect{b}$ and $\vect{c}
\in \IR^{d}$:
\begin{align}
  \vect{a} = \left[\vect{b}\right]\left[\vect{c}\right]
\end{align}
as $a_i = b_i c_i$ for all $i = 1,\dots d$ and for a vector~$\vect{a}\in
\IR^d$, we write for a set~$\possibleDOF$ with $|\possibleDOF| < d$:
\begin{align}
  \vect{a}_{\possibleDOF} = \vect{0}
\end{align}
as $a_p = 0$ for all $p = 1, \dots, |\possibleDOF|$, understood componentwise
respectively.

Moreover, we write for the scalar product for two scalar valued functions $f$,
$g \in L^2(\Omega_0)$
\begin{align}
  \left(f,g\right) = \int_{\Omega_0} f g \de \vect{X}_0,
\end{align}
for the scalar product for two vector valued functions $\vect{f}$, $\vect{g}
\in L^2(\Omega_0; \IR^d)$
\begin{align}
  \left(\vect{f}, \vect{g}\right) = \int_{\Omega_0} \vect{f} \cdot \vect{g} \de
  \vect{X}_0
\end{align}
and for the scalar product for two tensor valued functions $\tens{F}$,
$\tens{G} \in L^2(\Omega_0; \IR^{d,d})$
\begin{align}
  \left(\tens{F}, \tens{G}\right) = \int_{\Omega_0} \tens{F} \!:\!
  \tens{G} \de \vect{X}_0.
\end{align}
Boundary integrals for $\Gamma \subseteq \ptl\Omega_0$ are denoted with the
subscript of the respective boundary, e.g.,
\begin{align}
  \left(f,g\right)_{\Gamma} = \int_{\Gamma} f g \de \vect{S}_0.
\end{align}

%%%%%%%%%%%%%%%%%%%%%%%%%%%%%%%%%%%%%%%%%%%%%%%%%%%%%%%%%%%%%%%%%%%%%%%%%%%%%%%%
%% Appendix
%%%%%%%%%%%%%%%%%%%%%%%%%%%%%%%%%%%%%%%%%%%%%%%%%%%%%%%%%%%%%%%%%%%%%%%%%%%%%%%%

  %%%%%%%%%%%%%%%%%%%%%%%%%%%%%%%%%%%%%%%%%%%%%%%%%%%%%%%%%%%%%%%%%%%%%%%%%%%%%%
  %% Bibliography
  %%%%%%%%%%%%%%%%%%%%%%%%%%%%%%%%%%%%%%%%%%%%%%%%%%%%%%%%%%%%%%%%%%%%%%%%%%%%%%
  \bibliography{literature}

\begin{thebibliography}{10}
\expandafter\ifx\csname url\endcsname\relax
  \def\url#1{\texttt{#1}}\fi
\expandafter\ifx\csname urlprefix\endcsname\relax\def\urlprefix{URL }\fi
\expandafter\ifx\csname href\endcsname\relax
  \def\href#1#2{#2} \def\path#1{#1}\fi

\bibitem{tomaszewska2019lithium-ion}
A.~Tomaszewska, Z.~Chu, X.~Feng, S.~O'Kane, X.~Liu, J.~Chen, C.~Ji, E.~Endler,
  R.~Li, L.~Liu, Y.~Li, S.~Zheng, S.~Vetterlein, M.~Gao, J.~Du, M.~Parkes,
  M.~Ouyang, M.~Marinescu, G.~Offer, B.~Wu, Lithium-ion battery fast charging:
  {A} review, eTransportation 1 (2019) 100011.
\newblock \href {http://dx.doi.org/10.1016/j.etran.2019.100011}
  {\path{doi:10.1016/j.etran.2019.100011}}.

\bibitem{tian2015high}
H.~Tian, F.~Xin, X.~Wang, W.~He, W.~Han, High capacity group-{IV} elements
  ({Si}, {Ge}, {Sn}) based anodes for lithium-ion batteries, J. Materiomics
  1~(3) (2015) 153--169.
\newblock \href {http://dx.doi.org/10.1016/j.jmat.2015.06.002}
  {\path{doi:10.1016/j.jmat.2015.06.002}}.

\bibitem{li2021diverting}
P.~Li, H.~Kim, S.-T. Myung, Y.-K. Sun, Diverting exploration of silicon anode
  into practical way: A review focused on silicon-graphite composite for
  lithium ion batteries, Energy Stor. Mater. 35 (2021) 550--576.
\newblock \href {http://dx.doi.org/10.1016/j.ensm.2020.11.028}
  {\path{doi:10.1016/j.ensm.2020.11.028}}.

\bibitem{mo2020tin-graphene}
R.~Mo, X.~Tan, F.~Li, R.~Tao, J.~Xu, D.~Kong, Z.~Wang, B.~Xu, X.~Wang, C.~Wang,
  J.~Li, Y.~Peng, Y.~Lu, Tin-graphene tubes as anodes for lithium-ion batteries
  with high volumetric and gravimetric energy densities, Nat. Commun. 11~(1)
  (2020) 1374.
\newblock \href {http://dx.doi.org/10.1038/s41467-020-14859-z}
  {\path{doi:10.1038/s41467-020-14859-z}}.

\bibitem{zhang2011review}
W.-J. Zhang, A review of the electrochemical performance of alloy anodes for
  lithium-ion batteries, J. Power Sources 196~(1) (2011) 13--24.
\newblock \href {http://dx.doi.org/10.1016/j.jpowsour.2010.07.020}
  {\path{doi:10.1016/j.jpowsour.2010.07.020}}.

\bibitem{xu2016electrochemomechanics}
R.~Xu, K.~Zhao, Electrochemomechanics of electrodes in {Li}-ion batteries: {A}
  review, J. Electrochem. En. Conv. Stor. 13~(3) (2016) 030803.
\newblock \href {http://dx.doi.org/10.1115/1.4035310}
  {\path{doi:10.1115/1.4035310}}.

\bibitem{zhao2019review}
Y.~Zhao, P.~Stein, Y.~Bai, M.~Al-Siraj, Y.~Yang, B.-X. Xu, A review on modeling
  of electro-chemo-mechanics in lithium-ion batteries, J. Power Sources 413
  (2019) 259--283.
\newblock \href {http://dx.doi.org/10.1016/j.jpowsour.2018.12.011}
  {\path{doi:10.1016/j.jpowsour.2018.12.011}}.

\bibitem{song2015diffusion}
Y.~C. Song, Z.~Z. Li, A.~K. Soh, J.~Q. Zhang, Diffusion of lithium ions and
  diffusion-induced stresses in a phase separating electrode under
  galvanostatic and potentiostatic operations: Phase field simulations, Mech.
  Mater. 91 (2015) 363--371.
\newblock \href {http://dx.doi.org/10.1016/j.mechmat.2015.04.015}
  {\path{doi:10.1016/j.mechmat.2015.04.015}}.

\bibitem{delmas2008lithium}
C.~Delmas, M.~Maccario, L.~Croguennec, F.~Le~Cras, F.~Weill, Lithium
  deintercalation in {LiFePO$_\text{4}$} nanoparticles via a domino-cascade
  model, Nat. Mater. 7~(8) (2008) 665--671.
\newblock \href {http://dx.doi.org/10.1038/nmat2230}
  {\path{doi:10.1038/nmat2230}}.

\bibitem{van-der-ven2000phase}
A.~Van Der~Ven, C.~Marianetti, D.~Morgan, G.~Ceder, Phase transformations and
  volume changes in spinel {Li$_\text{x}$Mn$_\text{2}$O$_\text{4}$}, Solid
  State Ion. 135~(1--4) (2000) 21--32.
\newblock \href {http://dx.doi.org/10.1016/S0167-2738(00)00326-X}
  {\path{doi:10.1016/S0167-2738(00)00326-X}}.

\bibitem{walk2014comparison}
A.-C. Walk, M.~Huttin, M.~Kamlah, Comparison of a phase-field model for
  intercalation induced stresses in electrode particles of lithium ion
  batteries for small and finite deformation theory, Eur. J. Mech. A Solids 48
  (2014) 74--82.
\newblock \href {http://dx.doi.org/10.1016/j.euromechsol.2014.02.020}
  {\path{doi:10.1016/j.euromechsol.2014.02.020}}.

\bibitem{cahn1958free}
J.~W. Cahn, J.~E. Hilliard, Free energy of a nonuniform system. {I}.
  {Interfacial} free energy, J. Chem. Phys. 28~(2) (1958) 258--267.
\newblock \href {http://dx.doi.org/10.1063/1.1744102}
  {\path{doi:10.1063/1.1744102}}.

\bibitem{cahn1959free}
J.~W. Cahn, Free energy of a nonuniform system. {II}. {Thermodynamic} basis, J.
  Chem. Phys. 30~(5) (1959) 1121--1135.
\newblock \href {http://dx.doi.org/10.1063/1.1730145}
  {\path{doi:10.1063/1.1730145}}.

\bibitem{larche1973linear}
F.~Larch{\'{e}}, J.~W. Cahn, A linear theory of thermochemical equilibrium of
  solids under stress, Acta Metallurgica 21~(8) (1973) 1051--1063.
\newblock \href {http://dx.doi.org/10.1016/0001-6160(73)90021-7}
  {\path{doi:10.1016/0001-6160(73)90021-7}}.

\bibitem{garcke2001cahn-hilliard}
H.~Garcke, M.~Rumpf, U.~Weikard, The {C}ahn--{H}illiard equation with
  elasticity---finite element approximation and qualitative studies, Interfaces
  Free Bound. 3~(1) (2001) 101--118.
\newblock \href {http://dx.doi.org/10.4171/IFB/34} {\path{doi:10.4171/IFB/34}}.

\bibitem{garcke2005numerical}
H.~Garcke, U.~Weikard, Numerical approximation of the {C}ahn--{L}arch\'{e}
  equation, Numer. Math. 100~(4) (2005) 639--662.
\newblock \href {http://dx.doi.org/10.1007/s00211-004-0578-x}
  {\path{doi:10.1007/s00211-004-0578-x}}.

\bibitem{di-leo2014cahn-hilliard-type}
C.~V. Di~Leo, E.~Rejovitzky, L.~Anand, A {C}ahn--{H}illiard-type phase-field
  theory for species diffusion coupled with large elastic deformations:
  {Application} to phase-separating {Li}-ion electrode materials, J. Mech.
  Phys. Solids 70 (2014) 1--29.
\newblock \href {http://dx.doi.org/10.1016/j.jmps.2014.05.001}
  {\path{doi:10.1016/j.jmps.2014.05.001}}.

\bibitem{hennessy2020phase}
M.~G. Hennessy, A.~M\"{u}nch, B.~Wagner, Phase separation in swelling and
  deswelling hydrogels with a free boundary, Phys. Rev. E 101~(3) (2020)
  032501.
\newblock \href {http://dx.doi.org/10.1103/physreve.101.032501}
  {\path{doi:10.1103/physreve.101.032501}}.

\bibitem{werner2021multi-field}
M.~Werner, A.~Pandolfi, K.~Weinberg, A multi-field model for charging and
  discharging of lithium-ion battery electrodes, Contin. Mech. Thermodyn.
  33~(3) (2021) 661--685.
\newblock \href {http://dx.doi.org/10.1007/s00161-020-00943-8}
  {\path{doi:10.1007/s00161-020-00943-8}}.

\bibitem{huttin2012phase-field}
M.~Huttin, M.~Kamlah, Phase-field modeling of stress generation in electrode
  particles of lithium ion batteries, Appl. Phys. Lett. 101~(13) (2012)
  133902--1--133902--4.
\newblock \href {http://dx.doi.org/10.1063/1.4754705}
  {\path{doi:10.1063/1.4754705}}.

\bibitem{zhang2018nonlocal}
T.~Zhang, M.~Kamlah, A nonlocal species concentration theory for diffusion and
  phase changes in electrode particles of lithium ion batteries, Contin. Mech.
  Thermodyn. 30~(3) (2018) 553--572.
\newblock \href {http://dx.doi.org/10.1007/s00161-018-0624-z}
  {\path{doi:10.1007/s00161-018-0624-z}}.

\bibitem{castelli2021efficient}
G.~F. Castelli, L.~von Kolzenberg, B.~Horstmann, A.~Latz, W.~D{\"o}rfler,
  Efficient simulation of chemical-mechanical coupling in battery active
  particles, Energy Technol. 9~(6) (2021) 2000835.
\newblock \href {http://dx.doi.org/10.1002/ente.202000835}
  {\path{doi:10.1002/ente.202000835}}.

\bibitem{castelli2021numerical}
G.~F. Castelli, Numerical investigation of {C}ahn--{H}illiard-type phase-field
  models for battery active particles, Ph.D. thesis, Karlsruhe Institute of
  Technology (KIT) (2021).
\newblock \href {http://dx.doi.org/10.5445/IR/1000141249}
  {\path{doi:10.5445/IR/1000141249}}.

\bibitem{zhang2020mechanically}
T.~Zhang, M.~Kamlah, Mechanically coupled phase-field modeling of
  microstructure evolution in sodium ion batteries particles of
  {Na$_\text{x}$FePO$_\text{4}$}, J. Electrochem. Soc. 167~(2) (2020) 020508.
\newblock \href {http://dx.doi.org/10.1149/1945-7111/ab645a}
  {\path{doi:10.1149/1945-7111/ab645a}}.

\bibitem{wu2019phase}
L.~Wu, V.~De~Andrade, X.~Xiao, J.~Zhang, Phase field modeling of coupled phase
  separation and diffusion-induced stress in lithium iron phosphate particles
  reconstructed from synchrotron nano x-ray tomography, J. Electrochem. En.
  Conv. Stor. 16~(4) (2019) 041006.
\newblock \href {http://dx.doi.org/10.1115/1.4043155}
  {\path{doi:10.1115/1.4043155}}.

\bibitem{zhang2018sodium}
T.~Zhang, M.~Kamlah, Sodium ion batteries particles: {Phase}-field modeling
  with coupling of {C}ahn--{H}illiard equation and finite deformation
  elasticity, J. Electrochem. Soc. 165~(10) (2018) A1997--A2007.
\newblock \href {http://dx.doi.org/10.1149/2.0141810jes}
  {\path{doi:10.1149/2.0141810jes}}.

\bibitem{zhang2019phase-field}
T.~Zhang, M.~Kamlah, Phase-field modeling of the particle size and average
  concentration dependent miscibility gap in nanoparticles of
  {Li$_\text{x}$Mn$_\text{2}$O$_\text{4}$}, {Li$_\text{x}$FePO$_\text{4}$}, and
  {Na$_\text{x}$FePO$_\text{4}$} during insertion, Electrochim. Acta 298 (2019)
  31--42.
\newblock \href {http://dx.doi.org/10.1016/j.electacta.2018.12.007}
  {\path{doi:10.1016/j.electacta.2018.12.007}}.

\bibitem{chen2014phase-field}
L.~Chen, F.~Fan, L.~Hong, J.~Chen, Y.~Z. Ji, S.~L. Zhang, T.~Zhu, L.~Q. Chen, A
  phase-field model coupled with large elasto-plastic deformation:
  {Application} to lithiated silicon electrodes, J. Electrochem. Soc. 161~(11)
  (2014) F3164--F3172.
\newblock \href {http://dx.doi.org/10.1149/2.0171411jes}
  {\path{doi:10.1149/2.0171411jes}}.

\bibitem{zhang2019phase-field_1}
K.~Zhang, Y.~Li, F.~Wang, B.~Zheng, F.~Yang, A phase-field study of the effect
  of local deformation velocity on lithiation-induced stress in wire-like
  structures, J. Phys. D: Appl. Phys. 52 (2019) 145501.
\newblock \href {http://dx.doi.org/10.1088/1361-6463/ab00dc}
  {\path{doi:10.1088/1361-6463/ab00dc}}.

\bibitem{poluektov2018modelling}
M.~Poluektov, A.~B. Freidin, L.~Figiel, Modelling stress-affected chemical
  reactions in non-linear viscoelastic solids with application to lithiation
  reaction in spherical {S}i particles, Internat. J. Engrg. Sci. 128 (2018)
  44--62.
\newblock \href {http://dx.doi.org/10.1016/j.ijengsci.2018.03.007}
  {\path{doi:10.1016/j.ijengsci.2018.03.007}}.

\bibitem{kolzenberg2022chemo-mechanical}
L.~von Kolzenberg, A.~Latz, B.~Horstmann, Chemo-mechanical model of sei growth
  on silicon electrode particles, Batter. Supercaps 5~(2) (2022) e202100216.
\newblock \href {http://dx.doi.org/10.1002/batt.202100216}
  {\path{doi:10.1002/batt.202100216}}.

\bibitem{schoof2022parallelization}
R.~Schoof, G.~F. Castelli, W.~D\"orfler, Parallelization of a finite element
  solver for chemo-mechanical coupled anode and cathode particles in
  lithium-ion batteries, in: T.~Kvamsdal, K.~M. Mathisen, K.-A. Lie, M.~G.
  Larson (Eds.), 8th European Congress on Computational Methods in Applied
  Sciences and Engineering ({ECCOMAS} Congress 2022), {CIMNE}, 2022.
\newblock \href {http://dx.doi.org/10.23967/eccomas.2022.106}
  {\path{doi:10.23967/eccomas.2022.106}}.

\bibitem{laursen2002computational}
T.~A. Laursen, Computational contact and impact mechanics, Springer-Verlag,
  Berlin, Berlin, 2002.

\bibitem{wriggers2006computational}
P.~Wriggers, Computational contact mechanics, 2nd Edition, Springer, Berlin,
  2006.

\bibitem{willner2003kontinuums-}
K.~Willner, Kontinuums- und Kontaktmechanik: synthetische und analytische
  Darstellung, Engineering online library, Springer, Berlin, 2003.

\bibitem{alart1991mixed}
P.~Alart, A.~Curnier, A mixed formulation for frictional contact problems prone
  to {N}ewton like solution methods, Comput. Methods Appl. Mech. Engrg. 92~(3)
  (1991) 353--375.
\newblock \href {http://dx.doi.org/10.1016/0045-7825(91)90022-X}
  {\path{doi:10.1016/0045-7825(91)90022-X}}.

\bibitem{brunssen2007fast}
S.~Brunssen, F.~Schmid, M.~Sch\"{a}fer, B.~Wohlmuth, A fast and robust
  iterative solver for nonlinear contact problems using a primal-dual active
  set strategy and algebraic multigrid, Internat. J. Numer. Methods Engrg.
  69~(3) (2007) 524--543.
\newblock \href {http://dx.doi.org/10.1002/nme.1779}
  {\path{doi:10.1002/nme.1779}}.

\bibitem{fischer2005frictionless}
K.~A. Fischer, P.~Wriggers, Frictionless 2d contact formulations for finite
  deformations based on the mortar method, Comput. Mech. 36~(3) (2005)
  226--244.
\newblock \href {http://dx.doi.org/10.1007/s00466-005-0660-y}
  {\path{doi:10.1007/s00466-005-0660-y}}.

\bibitem{hintermuller2002primal-dual}
M.~Hinterm\"{u}ller, K.~Ito, K.~Kunisch, The primal-dual active set strategy as
  a semismooth {N}ewton method, SIAM J. Optim. 13~(3) (2002) 865--888 (2003).
\newblock \href {http://dx.doi.org/10.1137/S1052623401383558}
  {\path{doi:10.1137/S1052623401383558}}.

\bibitem{puso2004mortar}
M.~A. Puso, T.~A. Laursen, \href{https://doi.org/10.1016/j.cma.2003.10.010}{A
  mortar segment-to-segment contact method for large deformation solid
  mechanics}, Comput. Methods Appl. Mech. Engrg. 193~(6-8) (2004) 601--629.
\newblock \href {http://dx.doi.org/10.1016/j.cma.2003.10.010}
  {\path{doi:10.1016/j.cma.2003.10.010}}.
\newline\urlprefix\url{https://doi.org/10.1016/j.cma.2003.10.010}

\bibitem{wohlmuth2003monotone}
B.~I. Wohlmuth, R.~H. Krause, Monotone multigrid methods on nonmatching grids
  for nonlinear multibody contact problems, SIAM J. Sci. Comput. 25~(1) (2003)
  324--347.
\newblock \href {http://dx.doi.org/10.1137/S1064827502405318}
  {\path{doi:10.1137/S1064827502405318}}.

\bibitem{signorini1933sopra}
A.~Signorini, Sopra alcune questioni di elastostatica, Annali della Scuola
  Normale Superiore di Pisa - Scienze Fisiche e Matematiche 21~(2) (1933)
  143--148.

\bibitem{signorini1933sopra_1}
A.~Signorini, Sopra alcune questioni di statica dei sistemi continui, Annali
  della Scuola Normale Superiore di Pisa - Scienze Fisiche e Matematiche 2~(2)
  (1933) 231--251.

\bibitem{hueber2005primal-dual}
S.~H\"{u}eber, B.~I. Wohlmuth, A primal-dual active set strategy for non-linear
  multibody contact problems, Comput. Methods Appl. Mech. Engrg. 194~(27-29)
  (2005) 3147--3166.
\newblock \href {http://dx.doi.org/10.1016/j.cma.2004.08.006}
  {\path{doi:10.1016/j.cma.2004.08.006}}.

\bibitem{hueber2005priori}
S.~H\"{u}eber, M.~Mair, B.~I. Wohlmuth, A priori error estimates and an inexact
  primal-dual active set strategy for linear and quadratic finite elements
  applied to multibody contact problems, Appl. Numer. Math. 54~(3-4) (2005)
  555--576.
\newblock \href {http://dx.doi.org/10.1016/j.apnum.2004.09.019}
  {\path{doi:10.1016/j.apnum.2004.09.019}}.

\bibitem{hueber2013contact}
S.~H\"{u}eber, A.~Matei, B.~Wohlmuth, A contact problem for electro-elastic
  materials, ZAMM Z. Angew. Math. Mech. 93~(10-11) (2013) 789--800.
\newblock \href {http://dx.doi.org/10.1002/zamm.201200235}
  {\path{doi:10.1002/zamm.201200235}}.

\bibitem{hintermuller2003semismooth}
M.~Hinterm\"{u}ller, V.~A. Kovtunenko, K.~Kunisch, Semismooth newton methods
  for a class of unilaterally constrained variational problems, Technical
  Report 270{\;} Universit{\"a}t Graz/Technische Universit{\"a}t Graz. SFB
  F003-Optimierung und Kontrolle (2003).

\bibitem{frohne2016efficient}
J.~Frohne, T.~Heister, W.~Bangerth, Efficient numerical methods for the
  large-scale, parallel solution of elastoplastic contact problems, Internat.
  J. Numer. Methods Engrg. 105~(6) (2016) 416--439.
\newblock \href {http://dx.doi.org/10.1002/nme.4977}
  {\path{doi:10.1002/nme.4977}}.

\bibitem{hager2010semismooth}
C.~Hager, B.~I. Wohlmuth, Semismooth {N}ewton methods for variational problems
  with inequality constraints, GAMM-Mitt. 33~(1) (2010) 8--24.
\newblock \href {http://dx.doi.org/10.1002/gamm.201010002}
  {\path{doi:10.1002/gamm.201010002}}.

\bibitem{arndt2021deal-ii}
D.~Arndt, W.~Bangerth, B.~Blais, M.~Fehling, R.~Gassm{\"o}ller, T.~Heister,
  L.~Heltai, U.~K{\"o}cher, M.~Kronbichler, M.~Maier, P.~Munch, J.-P. Pelteret,
  S.~Proell, K.~Simon, B.~Turcksin, D.~Wells, J.~Zhang, The \texttt{deal.II}
  library, version 9.3, J. Numer. Math. 29~(3) (2021) 171--186.
\newblock \href {http://dx.doi.org/10.1515/jnma-2021-0081}
  {\path{doi:10.1515/jnma-2021-0081}}.

\bibitem{kornhuber2001adaptive}
R.~Kornhuber, R.~Krause, Adaptive multigrid methods for {S}ignorini's problem
  in linear elasticity, Comput. Vis. Sci. 4~(1) (2001) 9--20.
\newblock \href {http://dx.doi.org/10.1007/s007910100052}
  {\path{doi:10.1007/s007910100052}}.

\bibitem{sander2013towards}
O.~Sander, C.~Klapproth, J.~Youett, R.~Kornhuber, P.~Deuflhard, Towards an
  efficient numerical simulation of complex 3{D} knee joint motion, Comput.
  Vis. Sci. 16~(3) (2013) 119--138.
\newblock \href {http://dx.doi.org/10.1007/s00791-014-0227-6}
  {\path{doi:10.1007/s00791-014-0227-6}}.

\bibitem{de-los-reyes2012combined}
J.~C. De~Los~Reyes, S.~Gonz\'{a}lez~Andrade, A combined {BDF}-semismooth
  {N}ewton approach for time-dependent {B}ingham flow, Numer. Methods Partial
  Differ. Equ. 28~(3) (2012) 834--860.
\newblock \href {http://dx.doi.org/10.1002/num.20658}
  {\path{doi:10.1002/num.20658}}.

\bibitem{lauser2011new}
A.~Lauser, C.~Hager, R.~Helmig, B.~Wohlmuth, A new approach for phase
  transitions in miscible multi-phase flow in porous media, Adv. Water Resour.
  34~(8) (2011) 957--966.
\newblock \href {http://dx.doi.org/10.1016/j.advwatres.2011.04.021}
  {\path{doi:10.1016/j.advwatres.2011.04.021}}.

\bibitem{sa-ngiamsunthorn2021optimal}
P.~Sa~Ngiamsunthorn, A.~Suechoei, P.~Kumam, Optimal control for obstacle
  problems involving time-dependent variational inequalities with
  {L}iouville-{C}aputo fractional derivative, Adv. Differ. Equ. 2021 (2021)
  298.
\newblock \href {http://dx.doi.org/10.1186/s13662-021-03453-2}
  {\path{doi:10.1186/s13662-021-03453-2}}.

\bibitem{di-leo2015diffusion-deformation}
C.~V. Di~Leo, E.~Rejovitzky, L.~Anand, Diffusion-deformation theory for
  amorphous silicon anodes: The role of plastic deformation on electrochemical
  performance, Int. J. Solids Struct. 67-68 (2015) 283--296.
\newblock \href {http://dx.doi.org/10.1016/j.ijsolstr.2015.04.028}
  {\path{doi:10.1016/j.ijsolstr.2015.04.028}}.

\bibitem{castelli2021study}
G.~F. Castelli, W.~D{\"o}rfler, Study on an adaptive finite element solver for
  the {C}ahn--{H}illiard equation, in: F.~J. Vermolen, C.~Vuik (Eds.),
  Numerical Mathematics and Advanced Applications {ENUMATH} 2019, Vol. 139 of
  Lecture Notes in Computational Science and Engineering, Springer, Cham, 2021,
  pp. 245--253.
\newblock \href {http://dx.doi.org/10.1007/978-3-030-55874-1_23}
  {\path{doi:10.1007/978-3-030-55874-1_23}}.

\bibitem{holzapfel2000nonlinear}
G.~A. Holzapfel, Nonlinear Solid Mechanics, John Wiley \& Sons, Ltd.,
  Chichester, 2000.

\bibitem{braess2007finite}
D.~Braess, Finite Elements, 3rd Edition, Cambridge University Press, Cambridge,
  2007.
\newblock \href {http://dx.doi.org/10.1007/978-3-540-72450-6}
  {\path{doi:10.1007/978-3-540-72450-6}}.

\bibitem{latz2015multiscale}
A.~Latz, J.~Zausch, Multiscale modeling of lithium ion batteries: thermal
  aspects, Beilstein J. Nanotechnol. 6 (2015) 987--1007.
\newblock \href {http://dx.doi.org/10.3762/bjnano.6.102}
  {\path{doi:10.3762/bjnano.6.102}}.

\bibitem{latz2011thermodynamic}
A.~Latz, J.~Zausch, Thermodynamic consistent transport theory of {Li}-ion
  batteries, J. Power Sources 196~(6) (2011) 3296--3302.
\newblock \href {http://dx.doi.org/10.1016/j.jpowsour.2010.11.088}
  {\path{doi:10.1016/j.jpowsour.2010.11.088}}.

\bibitem{schammer2021theory}
M.~Schammer, B.~Horstmann, A.~Latz, Theory of transport in highly concentrated
  electrolytes, J. Electrochem. Soc. 168~(2) (2021) 026511.
\newblock \href {http://dx.doi.org/10.1149/1945-7111/abdddf}
  {\path{doi:10.1149/1945-7111/abdddf}}.

\bibitem{anand2012cahn-hilliard-type}
L.~Anand, A {C}ahn--{H}illiard-type theory for species diffusion coupled with
  large elastic-plastic deformations, J. Mech. Phys. Solids 60~(12) (2012)
  1983--2002.
\newblock \href {http://dx.doi.org/10.1016/j.jmps.2012.08.001}
  {\path{doi:10.1016/j.jmps.2012.08.001}}.

\bibitem{zhang2018lithiation-induced}
K.~Zhang, Y.~Li, J.~Wu, B.~Zheng, F.~Yang, Lithiation-induced buckling of
  wire-based electrodes in lithium-ion batteries: A phase-field model coupled
  with large deformation, Int. J. Solids Struct. 144-145 (2018) 289--300.
\newblock \href {http://dx.doi.org/10.1016/j.ijsolstr.2018.05.014}
  {\path{doi:10.1016/j.ijsolstr.2018.05.014}}.

\bibitem{chan2007high-performance}
C.~K. Chan, H.~Peng, G.~Liu, K.~McIlwrath, X.~F. Zhang, R.~A. Huggins, Y.~Cui,
  High-performance lithium battery anodes using silicon nanowires, Nat.
  Nanotechnol. 3~(1) (2007) 31--35.
\newblock \href {http://dx.doi.org/10.1038/nnano.2007.411}
  {\path{doi:10.1038/nnano.2007.411}}.

\bibitem{keil2016calendar}
P.~Keil, S.~F. Schuster, J.~Wilhelm, J.~Travi, A.~Hauser, R.~C. Karl,
  A.~Jossen, Calendar aging of lithium-ion batteries, J. Electrochem. Soc.
  163~(9) (2016) A1872--A1880.
\newblock \href {http://dx.doi.org/10.1149/2.0411609jes}
  {\path{doi:10.1149/2.0411609jes}}.

\bibitem{latz2013thermodynamic}
A.~Latz, J.~Zausch, Thermodynamic derivation of a {Butler--Volmer} model for
  intercalation in {Li}-ion batteries, Electrochim. Acta 110 (2013) 358--362.
\newblock \href {http://dx.doi.org/10.1016/j.electacta.2013.06.043}
  {\path{doi:10.1016/j.electacta.2013.06.043}}.

\bibitem{hoffmann2018influence}
V.~Hoffmann, G.~Pulletikurthi, T.~Carstens, A.~Lahiri, A.~Borodin, M.~Schammer,
  B.~Horstmann, A.~Latz, F.~Endres, Influence of a silver salt on the
  nanostructure of a {Au}(111)/ionic liquid interface: {An} atomic force
  microscopy study and theoretical concepts, Phys. Chem. Chem. Phys. 20~(7)
  (2018) 4760--4771.
\newblock \href {http://dx.doi.org/10.1039/C7CP08243F}
  {\path{doi:10.1039/C7CP08243F}}.

\bibitem{friedman1982variational}
A.~Friedman, Variational principles and free-boundary problems, Pure and
  applied mathematics, Wiley, New York, 1982.

\bibitem{kornhuber1997adaptive}
R.~Kornhuber, Adaptive Monotone Multigrid Methods for Nonlinear Variational
  Problems, 1st Edition, Advances in Numerical Mathematics, B. G. Teubner,
  Stuttgart, 1997.

\bibitem{haslinger1980contact}
J.~Haslinger, I.~Hlav{\'{a}}{\v{c}}ek, Contact between elastic bodies. {I}.
  {C}ontinuous problems, Appl. Math. 25~(5) (1980) 324--347.
\newblock \href {http://dx.doi.org/10.21136/am.1980.103868}
  {\path{doi:10.21136/am.1980.103868}}.

\bibitem{boieri1987existence}
P.~Boieri, F.~Gastaldi, D.~Kinderlehrer, Existence, uniqueness, and regularity
  results for the two-body contact problem, Appl. Math. Optim. 15~(3) (1987)
  251--277.
\newblock \href {http://dx.doi.org/10.1007/BF01442654}
  {\path{doi:10.1007/BF01442654}}.

\bibitem{hlavacek1988solution}
I.~Hlav\'{a}\v{c}ek, J.~Haslinger, J.~Ne\v{c}as, J.~Lov\'{\i}\v{s}ek, Solution
  of variational inequalities in mechanics, Vol.~66 of Applied Mathematical
  Sciences, Springer-Verlag, New York, 1988.
\newblock \href {http://dx.doi.org/10.1007/978-1-4612-1048-1}
  {\path{doi:10.1007/978-1-4612-1048-1}}.

\bibitem{kinderlehrer2000introduction}
D.~Kinderlehrer, G.~Stampacchia, An introduction to variational inequalities
  and their applications, Vol.~31 of Classics in Applied Mathematics, Society
  for Industrial and Applied Mathematics (SIAM), Philadelphia, PA, 2000.
\newblock \href {http://dx.doi.org/10.1137/1.9780898719451}
  {\path{doi:10.1137/1.9780898719451}}.

\bibitem{ben-belgacem1999extension}
F.~Ben~Belgacem, P.~Hild, P.~Laborde, Extension of the mortar finite element
  method to a variational inequality modeling unilateral contact, Math. Models
  Methods Appl. Sci. 9~(2) (1999) 287--303.
\newblock \href {http://dx.doi.org/10.1142/S0218202599000154}
  {\path{doi:10.1142/S0218202599000154}}.

\bibitem{hild2000numerical}
P.~Hild, Numerical implementation of two nonconforming finite element methods
  for unilateral contact, Comput. Methods Appl. Mech. Engrg. 184~(1) (2000)
  99--123.
\newblock \href {http://dx.doi.org/10.1016/S0045-7825(99)00096-1}
  {\path{doi:10.1016/S0045-7825(99)00096-1}}.

\bibitem{wohlmuth2000mortar}
B.~I. Wohlmuth, A mortar finite element method using dual spaces for the
  {L}agrange multiplier, SIAM J. Numer. Anal. 38~(3) (2000) 989--1012.
\newblock \href {http://dx.doi.org/10.1137/S0036142999350929}
  {\path{doi:10.1137/S0036142999350929}}.

\bibitem{reichelt1997matlab}
M.~W. Reichelt, L.~F. Shampine, J.~Kierzenka,
  \href{http://www.mathworks.com}{Matlab \texttt{ode15s}}, copyright 1984--2020
  {The MathWorks, Inc.} (1997).
\newline\urlprefix\url{http://www.mathworks.com}

\bibitem{shampine1997matlab}
L.~F. Shampine, M.~W. Reichelt, The {MATLAB} {ODE} suite, SIAM J. Sci. Comput.
  18~(1) (1997) 1--22.
\newblock \href {http://dx.doi.org/10.1137/S1064827594276424}
  {\path{doi:10.1137/S1064827594276424}}.

\bibitem{shampine1999solving}
L.~F. Shampine, M.~W. Reichelt, J.~A. Kierzenka, Solving index-{$1$} {DAE}s in
  {MATLAB} and {S}imulink, SIAM Rev. 41~(3) (1999) 538--552.
\newblock \href {http://dx.doi.org/10.1137/S003614459933425X}
  {\path{doi:10.1137/S003614459933425X}}.

\bibitem{shampine2003solving}
L.~F. Shampine, I.~Gladwell, S.~Thompson, Solving {ODE}s with {MATLAB},
  Cambridge University Press, Cambridge, 2003.
\newblock \href {http://dx.doi.org/10.1017/CBO9780511615542}
  {\path{doi:10.1017/CBO9780511615542}}.

\bibitem{ainsworth2000posteriori}
M.~Ainsworth, J.~T. Oden, A Posteriori Error Estimation in Finite Element
  Analysis, Pure and Applied Mathematics, John Wiley \& Sons, Inc., New York,
  2000.

\bibitem{banas2008adaptive}
L.~Ba\v{n}as, R.~N{\"u}rnberg, Adaptive finite element methods for
  {C}ahn--{H}illiard equations, J. Comput. Appl. Math. 218~(1) (2008) 2--11.
\newblock \href {http://dx.doi.org/10.1016/j.cam.2007.04.030}
  {\path{doi:10.1016/j.cam.2007.04.030}}.

\bibitem{team2020trilinos}
T.~{T}rilinos~{P}roject {T}eam, \href{https://trilinos.github.io}{The
  {T}rilinos {P}roject {W}ebsite} (2020).
\newline\urlprefix\url{https://trilinos.github.io}

\bibitem{davis2004algorithm}
T.~A. Davis, Algorithm 832: {UMFPACK} {V}4.3---an unsymmetric-pattern
  multifrontal method, ACM Trans. Math. Software 30~(2) (2004) 196--199.
\newblock \href {http://dx.doi.org/10.1145/992200.992206}
  {\path{doi:10.1145/992200.992206}}.

\bibitem{al-obeidi2016mechanical}
A.~Al-Obeidi, D.~Kramer, S.~T. Boles, R.~M\"{o}nig, C.~V. Thompson, Mechanical
  measurements on lithium phosphorous oxynitride coated silicon thin film
  electrodes for lithium-ion batteries during lithiation and delithiation,
  Appl. Phys. Lett. 109~(7) (2016) 071902.
\newblock \href {http://dx.doi.org/10.1063/1.4961234}
  {\path{doi:10.1063/1.4961234}}.

\end{thebibliography}

\end{document}